%&Plain
%%%%%%%%%%%%%%%%%%%%%%%%%%%%%%%%%%%%%%%%%%%%%%%%%%%%%%%%%%%%%%%%%%%%%%%%%%
%%%%%%%%%%%%%%% contatori fino a quattro entrate
%%%%%%%%%%%%%%% uso: \N o \NN o \NNN o \NNNN
%%%%%%%%%%%%%%% effetto: stampa numerazione nel font corrente
%%%%%%%%%%%%%%%%%%%%%%%%%%%%%%%%%%%%%%%%%%%%%%%%%%%%%%%%%%%%%%%%%%%%%%%%%%

\newcount\s
\newcount\n
\def\clearn{\n =0}
\def\N{\global\advance\n by 1 \global\s=1
                       \global\clearnn
                       \global\clearnnn
                       \global\clearnnnn
                       {\the\n}%
                     }
\newcount\nn
\def\clearnn{\nn =0}
\def\NN{\global\advance\nn by 1 \global\s=2
                        \global\clearnnn
                        \global\clearnnnn
                        {\the\n}.{\the\nn}%
                   }
\newcount\nnn
\def\clearnnn{\nnn =0}
\def\NNN{\global\advance\nnn by 1 \global\s=3
                        \global\clearnnnn
                        {\the\n}.{\the\nn}.{\the\nnn}%
                     }

\newcount\nnnn
\def\clearnnnn{\nnnn =0}
\def\NNNN{\global\advance\nnnn by 1 \global\s=4
                        {\the\n}.{\the\nn}.{\the\nnn}.{\the\nnnn}%
                     }

%%%%%%%%%%%%%%%%%%%%%%%%%%%%%%%%%%%%%%%%%%%%%%%%%%%%%%%%%%%%%%%%%%%%%
%%%%%%%%%%%%%% \rif{NOME} crea le label per il riferimento
%%%%%%%%%%%%%%  (lo stato della numerazione in quel momento)
%%%%%%%%%%%%%% \cite{NOME} scrive la citazione corrispondente a NOME
%%%%%%%%%%%%%%  (nel font corrente)
%%%%%%%%%%%%%%%%%%%%%%%%%%%%%%%%%%%%%%%%%%%%%%%%%%%%%%%%%%%%%%%%%%%%%

\def\rif#1{\expandafter\xdef\csname s#1\endcsname{\number\s}%
            \ifnum\s>0 \expandafter\xdef\csname n#1\endcsname{\number\n}\fi%
            \ifnum\s>1 \expandafter\xdef\csname nn#1\endcsname{\number\nn}\fi%
            \ifnum\s>2 \expandafter\xdef\csname nnn#1\endcsname{\number\nnn}\fi%
            \ifnum\s>3 \expandafter\xdef\csname nnnn#1\endcsname{\number\nnnn}\fi%
                   }

\def\cite#1{%
            \expandafter\ifcase\csname s#1\endcsname%
            \or \csname n#1\endcsname%
            \or \csname n#1\endcsname.\csname nn#1\endcsname%
            \or \csname n#1\endcsname.\csname nn#1\endcsname.\csname nnn#1\endcsname%
            \or \csname n#1\endcsname.\csname nn#1\endcsname.\csname nnn#1\endcsname.%
                \csname nnnn#1\endcsname%
                   \fi}

%%%%%%%%%%%%%%%%%%%%%%%%%%%%%%%%%%%%%%%%%%%%%%%%%%%%%%%%%%%%%%%%%%%%%%%%%
\font \proc=cmcsc10
\def\proclaim#1. #2\par
                  {\medbreak  {\proc #1. \enspace }{\sl #2\par }
                    \ifdim \lastskip <\medskipamount \removelastskip
                            \penalty 55\medskip \fi}
%%%%%%%%%%%%%%%%%%%%%%%%%%%%%%%%%%%%%%%%%%%%%%%%%%%%%%%%%%%%%%%%%%%%%%%%%

%%%%%%%%%%%%%%%%%%%%%%%%%%%%%%%%%%%%%%%%%%%%%%%%%%%%%%%%%%%%%%%%%%%%%%%%%
\def\biblio#1#2{ \item{\hbox to 1.2cm{[#1]\hss}}
                           {#2} }
%%%%%%%%%%%%%%%%%%%%%%%%%%%%%%%%%%%%%%%%%%%%%%%%%%%%%%%%%%%%%%%%%%%%%%%%%

%%%%%%%%%%%%%%%%%%%%%%%%%%%%%%%%%%%%%%%%%%%%%%%%%%%%%%%%%%%%%%%%%%%%%%%%%
%% ambienti predefiniti #1=numerazione (N, NN, NNN o NNNN)
%%                      #2=testo

\def\numfont{\bf}
%\font\numfont=?

\def\prop#1#2{\proclaim{{\numfont\csname#1\endcsname.} Proposition}. {#2}\par}
\def\theo#1#2{\proclaim{{\numfont\csname#1\endcsname.} Theorem}. {#2}\par}
\def\lemm#1#2{\proclaim{{\numfont\csname#1\endcsname.} Lemma}. {#2}\par}
\def\defi#1#2{\proclaim{{\numfont\csname#1\endcsname.} Definition}. {#2}\par}
\def\coro#1#2{\proclaim{{\numfont\csname#1\endcsname.} Corollary}. {#2}\par}

\def\Prop#1#2#3{% 
\proclaim{{\numfont\csname#1\endcsname.} Proposition (#2)}. {#3}\par}
\def\Theo#1#2#3{%
\proclaim{{\numfont\csname#1\endcsname.} Theorem (#2)}. {#3}\par}
\def\Lemm#1#2#3{%
\proclaim{{\numfont\csname#1\endcsname.} Lemma (#2)}. {#3}\par}
\def\Defi#1#2#3{%
\proclaim{{\numfont\csname#1\endcsname.} Definition (#2)}. {#3}\par}
\def\Coro#1#2#3{% 
\proclaim{{\numfont\csname#1\endcsname.} Corollary (#2)}. {#3}\par}

\long\def\proo#1#2{{{\numfont\csname#1\endcsname.}\proc \ Proof.}
                                    #2
                                    \hfill$\square$\vskip 10pt
                                }
\long\def\exam#1#2{{{\numfont\csname#1\endcsname.}\proc \ Example.}
                                    #2
                                    \vskip 10pt
                                }

\long\def\rema#1#2{{{\numfont\csname#1\endcsname.}\proc \ Remark.}
                                    #2
                                    \vskip 10pt
                                }

\long\def\Rema#1#2#3{{{\numfont\csname#1\endcsname.}\proc \ Remark (#2).}
                                    #3
                                    \vskip 10pt
                                }

\long\def\nota#1#2{{{\numfont\csname#1\endcsname.}\proc \ Notations.}
                                    #2
                                    \vskip 10pt
                                }

%%%%%%%%%%%%%%%%%%%%%%%%%%%%%%%%%%%%%%%%%
%% ambienti predefiniti #1=titolo

\font\titlefont=cmbx10 scaled \magstep2
\font\authorfont=cmcsc10 scaled \magstep1
\font\sectionfont=cmbx10 scaled \magstep1
\font\subsectionfont=cmbx10
\font\abstractfont=cmr8
\font\titleabstractfont=cmcsc8

\def\section#1{\vskip 15pt
                          {\sectionfont\N.\ #1}
                          \vskip 10pt}
\def\subsection#1{\vskip 12pt
                             {\subsectionfont\NN.\ #1}
                             \vskip 8pt}
\def\intro{\vskip 15pt
                          {\sectionfont Introduction.}
                          \vskip 10pt}
\def\refer{\vskip 15pt
                          {\sectionfont References.}
                          \vskip 10pt}
\def\title#1{\vskip 5pt
                          \centerline{\titlefont #1}
                          \vskip 10pt}
\def\author#1{\vskip 10pt
                          %\centerline{by}\vskip 5pt
                          \centerline{\authorfont #1}
                          \vskip 20pt}

\def\abstract#1{\vskip 10pt
\centerline{\vtop{\hsize=12truecm\baselineskip=9pt\strut
                 {\titleabstractfont Abstract. }
                         \abstractfont #1}}
                          \vskip 20pt}

%%%%%%%%%%%%%%%%%%%%%%%%%%%%% start %%%%%%%%%%%%%%%%%%%%%%%%%%%%%%
\clearn

%%%%%%%%%%%%%%%%%%%%%%%%%%%% DIMENSIONS
\vsize= 230truemm
\hsize= 165truemm
\voffset= 0truemm
\hoffset= -3truemm
\topskip=15pt 

\abovedisplayskip=6pt plus2pt minus 4pt
\belowdisplayskip=6pt plus2pt minus 4pt
\abovedisplayshortskip=0pt plus2pt
\belowdisplayshortskip=2pt plus1pt minus 1pt

%%%%%%%%%%%%%%%%%%%%%%%%%%%

%%%%%%%%%%%%%%%%%%%%%%%%%%%%%%%%%%%%%%%%%%%%%%%%%%%%%%%%%%
\font \small=cmr8

\font\smalltt=cmtt8

%%%%%%%%%%%%%%%%%%%%%%%%%%%%%%%%%%%%%%%%%%%%%%%%%%%%%%%%%%

%%%%%%%%%%%%%%% Generic by Maurizio %%%%%%%%%%%%%%%%%%%%%%%%%%%%%%

%% eliminare gli overfull

\def\NoBlackBoxes{\global\overfullrule 0pt}
\NoBlackBoxes

%% Abbreviazioni

\def\sss{\scriptscriptstyle}

\def\ov{\overline}

\def\r{\right}
\def\l{\left}

\def\phi{\varphi}

\def\theta{\vartheta}

\def\rho{\varrho}

%% Font

%\font \title=cmbx10 scaled \magstep5
%\font \tit=cmbx10 scaled \magstep2
%\font \tpar=cmbx10 scaled \magstep1

%% font matematici: bold, gotic, calligrafico

%%%%%%%%%%%%%%%%%%%%%%%%%%%%%%%%%%%%%%%%%%%%%%%%%%%%%%%%%%%%%%%%%%%%%%%%%%
% The following allows the use of Ralph Smith's Formal Script symbols
% in Plain TeX documents.  Use \scr like \cal.
% Set the font sizes and restore the `at' clauses if you want them bigger.
% You can use this method in LaTeX, but only at one basic size.
% If you need symbols in LaTeX titles, captions, etc., work it out or ask
% a LaTeXpert.

\font\tenscr=rsfs10 % scaled \magstep1
\font\sevenscr=rsfs7 % scaled \magstep1
\font\fivescr=rsfs5 % scaled \magstep1
\skewchar\tenscr='177 \skewchar\sevenscr='177 \skewchar\fivescr='177
\newfam\scrfam \textfont\scrfam=\tenscr \scriptfont\scrfam=\sevenscr
\scriptscriptfont\scrfam=\fivescr
\def\scr{\fam\scrfam}
%%%%%%%%%%%%%%%%%%%%%%%%%%%%%%%%%%%%%%%%%%%%%%%%%%%%%%%%%%%%%%%%%%%%%%%%%%%

\def\b#1{ {\Bbb {#1}} }

\def\c#1{ {\cal #1} }
\def\c#1{ {\scr #1} }

%% Limiti
%% vedere comunque quelli AMS
%% la prima serie richiede UnOvAr.tex

\def\limpr{\mathop{\underleftarrow{\rm lim}}}
\def\uplimpr{\mathop{\hbox{``}\underleftarrow{\rm lim}\hbox{''}}}
\def\limin{\mathop{\underrightarrow{\rm lim}}}
\def\uplimin{\mathop{\hbox{``}\underrightarrow{\rm lim}\hbox{''}}}

%% Nomi

\def\Hom{{\rm Hom}}

\def\h{{\rm h}}
\def\ob{ {\rm ob} }

\def\Ind{{\rm Ind}}
\def\Pro{{\rm Pro}}
\def\ind{{\rm ind}}
\def\pro{{\rm pro}}
\def\opp{{o}}

\def\bC{{\bf C}}
\def\bK{{\bf K}}
\def\bD{{\bf D}}

%% Trucco: come incastellare 

%%% File: Laps.TeX
\catcode`\@=11   % Let's pretend @ is a letter 
%         Vertical `laps'; cf. \llap and \rlap
\long\def\ulap#1{\vbox to \z@{\vss#1}}
\long\def\dlap#1{\vbox to \z@{#1\vss}}
     
%         And centered horizontal and vertical `laps'
\def\xlap#1{\hbox to \z@{\hss#1\hss}}
\long\def\ylap#1{\vbox to \z@{\vss#1\vss}}
     
%         And a `lap' centered on its midpoint

%         And the not-long definitions
\def\ulap#1{\vbox to \z@{\vss#1}}
\def\dlap#1{\vbox to \z@{#1\vss}}
\def\hlap#1{\hbox to \z@{\hss#1\hss}}
\def\vlap#1{\vbox to \z@{\vss\hbox{#1}\vss}}
\def\clap#1{\vbox to \z@{\vss\hbox to \z@{\hss#1\hss}\vss}}

\catcode`\@=12   % That's enough pretending for one day.

\catcode`\@=11  

\def\mmatrix#1{\null\,\vcenter{\normalbaselines\m@th 
               \ialign{\hfil$##$\hfil&&\ \hfil$##$\hfil\crcr 
               \mathstrut\crcr\noalign{\kern-\baselineskip } 
               #1\crcr\mathstrut\crcr\noalign{\kern-\baselineskip}}}\,}

\catcode`\@=12

\def\diagram{%
    \def\normalbaselines{\baselineskip15pt\lineskip3pt\lineskiplimit3pt}%
    \mmatrix
}

%%%%%%%%%%%%%%%%%%%%%%%%%%%%%%%%%%%%%%%

\input amssym.def
\input amssym.tex

%%%%%%%%%%%%%%%%%%%%%%%%%%%%%%%%%%%%%%%

%%%%%%%%%%%%%%%%%%%%%%%%%%% TeX code
\catcode`\@=11 
%% see \rightarrowfill, \leftarrowfill, 
%%     \overrightarrow, \overleftarrow 
%% of TeX

\def\Rightarrowfill{ $\m@th \mathord =\mkern -6mu\cleaders 
  \hbox {$\mkern -2mu\mathord =\mkern -2mu$}
     \hfill \mkern -6mu\mathord \Rightarrow $}
\def\Leftarrowfill{$\m@th \mathord \Leftarrow \mkern -6mu\cleaders 
        \hbox {$\mkern -2mu\mathord =\mkern -2mu$}
        \hfill \mkern -6mu\mathord =$}
\def\Leftrightarrowfill{$\m@th \mathord \Leftarrow \mkern -6mu\cleaders 
        \hbox {$\mkern -2mu\mathord =\mkern -2mu$}
        \hfill \mkern -6mu\mathord \Rightarrow $}

\def\leftrightarrowfill{$\m@th \mathord \leftarrow \mkern -6mu\cleaders 
        \hbox {$\mkern -2mu\mathord -\mkern -2mu$}
        \hfill \mkern -6mu\mathord \rightarrow $}

\def\underrightarrow#1{\vtop {\m@th \ialign {##\crcr 
                $\hfil \displaystyle { #1 }\hfil $\crcr 
                \noalign {\kern 1pt \nointerlineskip } 
               \rightarrowfill \crcr \noalign {\kern 0pt }}}}

\def\underleftarrow#1{\vtop {\m@th \ialign {##\crcr 
                $\hfil \displaystyle { #1 }\hfil $\crcr 
                \noalign {\kern 1pt \nointerlineskip } 
               \leftarrowfill \crcr \noalign {\kern 0pt }}}}

\def\underleftrightarrow#1{\vtop {\m@th \ialign {##\crcr 
                $\hfil \displaystyle { #1 }\hfil $\crcr 
                \noalign {\kern 0pt \nointerlineskip } 
               \leftrightarrowfill \crcr \noalign {\kern 0pt }}}}

\def\overleftrightarrow#1{\vbox {\m@th \ialign {##\crcr 
              \noalign {\kern 3pt } \leftrightarrowfill \crcr 
              \noalign {\kern 3pt \nointerlineskip } 
              $\hfil \displaystyle { #1 }\hfil$\crcr }}}

%% see \upbracefill, \downbracefill, 
%%     \underbrace, \overbrace 
%% of TeX

\def\uprotfill{$\m@th \bracelu \leaders \vrule \hfill \braceru $}
\def\downrotfill{$\m@th \braceld \leaders \vrule \hfill \bracerd $}
\def\underrot#1{\vtop {\m@th \ialign {##\crcr 
                $\hfil \displaystyle {\ #1 \ }\hfil $\crcr 
                \noalign {\kern 3pt \nointerlineskip } 
               \uprotfill \crcr \noalign {\kern 3pt }}}}
\def\overrot#1{\vbox {\m@th \ialign {##\crcr 
              \noalign {\kern 3pt } \downrotfill \crcr 
              \noalign {\kern 3pt \nointerlineskip } 
              $\hfil \displaystyle {\ #1 \ }\hfil$\crcr }}}

%% see \hrulefill, 
%%     \underline, \overline 
%% of TeX

\def\trattfill#1#2{\leaders \hbox{\hbox to#1pt{\hrulefill}
                                        \hskip #2pt} \hfill }
%% #1 width of line
%% #2 width of space (also negative) 

\def\undertratt#1#2#3{\vtop {\ialign {##\crcr 
                \hfil#3\hskip 6pt \hfil \crcr 
                \noalign {\kern 3pt \nointerlineskip } 
                \trattfill{#1}{#2} \crcr \noalign {\kern 3pt }}}}
\def\overtratt#1#2#3{\vbox {\ialign {##\crcr 
              \noalign {\kern 3pt } \trattfill{#1}{#2} \crcr 
              \noalign {\kern 3pt \nointerlineskip } 
              \hfil#3\hskip 6pt \hfil\crcr }}}

\def\dotfill#1{\leaders \hbox to#1pt{\hss . \hss} \hfill }
\def\cdotfill#1{\leaders \hbox to#1pt{\hss $\cdot$ \hss} \hfill }

%% #1 width of space

%% find on the net: undertilde

\def\undertilde#1{\mathord{\vtop{\ialign{##\crcr
   $\hfil\displaystyle{#1}\hfil$\crcr\noalign{\kern1.5pt\nointerlineskip}
   $\hfil\tilde{}\hfil$\crcr\noalign{\kern1.5pt}}}}}

\catcode`\@=12
%%%%%%%%%%%%%%%%%%%%%%%%%%%%%%%%%%%%%%%%%

%%%%%%%%%%%%%%%%%%%%%%%%%%%%%%%%%%%%%%%
%%%%%%% Lenght Variable Arrows 
%%%%%%% this for the construction

\def\rightto#1{\hbox to#1pt{\rightarrowfill}}
\def\leftto#1{\hbox to#1pt{\leftarrowfill}}
\def\leftrightto#1{\hbox to#1pt{\leftrightarrowfill}}

\def\mapstoto#1{\mapstochar\mkern -4mu\hbox to#1pt{\rightarrowfill}}
\def\rmapstoto#1{\hbox to#1pt{\leftarrowfill}\mkern -4mu\mapstochar}

\def\hookrightto#1{\lhook\mkern -8mu\hbox to#1pt{\rightarrowfill}}
\def\hookleftto#1{\hbox to#1pt{\leftarrowfill}\mkern -8mu\rhook}

\def\twoheadrightto#1{\hbox to#1pt{\rightarrowfill}\mkern -20mu\rightarrow }
\def\twoheadleftto#1{\leftarrow\mkern -20mu\hbox to#1pt{\leftarrowfill}}

\def\hooktwoheadrightto#1{\lhook\mkern -8mu%
                          \hbox to#1pt{\rightarrowfill}%
                          \mkern -20mu\rightarrow }
\def\hooktwoheadleftto#1{\leftarrow\mkern -20mu%
                         \hbox to#1pt{\leftarrowfill}%
                         \mkern -8mu\rhook}

%%%%%%% user Arrows 
%%
%% these arrows again call for the lenght in pt
%% R=right L=left     
%% h=hook  t=twohead   to=to

\def\R#1{\mathop{\rightto{#1}}\limits}
\def\L#1{\mathop{\leftto{#1}}\limits}

\def\U#1{\left\uparrow\vbox to#1pt{}\right.}
\def\D#1{\left\downarrow\vbox to#1pt{}\right.}

%%
%% here w=big-width W=Big ww=bigg WW=Bigg

\def\w{20.5}   \def\h{7.5}
\def\W{25.5}   \def\H{12.5}
\def\ww{30.5}  
\def\WW{35.5}  

%% to use as in $\R\W_{}^{}$ 

%%%%%% to end, how to put one arrow over another? 

\def\drla#1{\mathop{
              \raise 3pt\hbox to0pt{\rightto{#1}\hss}%
              \lower 3pt\leftto{#1}}\limits}

\def\dRL#1{\mathop{
              \raise 3pt\hbox to0pt{\rightto{#1}\hss}%
              \lower 3pt\leftto{#1}}\limits}

\def\dra#1{\mathop{
              \raise 3pt\hbox to0pt{\rightto{#1}\hss}%
              \lower 3pt\rightto{#1}}\limits}

\def\RR#1{\mathop{
              \raise 3pt\hbox to0pt{\rightto{#1}\hss}%
              \lower 3pt\rightto{#1}}\limits}
\def\UU#1{\U{#1}\kern-5pt\U{#1}}

%% it's enough, for the moment!

%%%%%%%%%%%%%%%%% local definitions

%%%%%%%%%%%%%%%%%%%%%%%%%%%% arrows and names 
\def\ra{\mathop{\rightarrow}\limits}
\def\la{\mathop{\leftarrow}\limits}

\def\lra{\mathop{\longrightarrow}\limits}

\def\nea{\mathop{\nearrow}\limits}
\def\sea{\mathop{\searrow}\limits}

\def\ln#1{\llap{$\scriptstyle #1$}}  %left name
\def\rn#1{\rlap{$\scriptstyle #1$}}  %right name

\def\hw{\hidewidth}

%%%%%%%%%%%%%%%%%%%%%%%%%% categories 
\def\Set{{{\cal S}et}}
\def\Ab{{{\cal A}b}}

\def\Funct{{\c F\! unct}}

\def\cov{ {\sss \vee} }

\def\con{ {\sss \wedge} }

%%%%%%%%%%%%%%%%%%%%%%%%%%%%%%% symbols

\def\id{{\rm id}}
\def\bydef{\mathrel{:=}}

\def\isom{\cong}

\def\int{{\rm int}}

\def\epsilon{\varepsilon}

\def\point{{\bullet}}

%%%%%%%%%%%%%%%%%%%%%%%%%%%%%%%%%%%%

\abovedisplayskip=6pt plus2pt minus 4pt 
\belowdisplayskip=6pt plus2pt minus 4pt 
\abovedisplayshortskip=0pt plus2pt  
\belowdisplayshortskip=2pt plus1pt minus 1pt

%%%%%%%%%%%%%%%%%%%%%%%%%%%%%%%%%%%%%%%%
\font \small=cmr8
\font\smalltt=cmtt8

\headline{ 
%\hss\small 
%DLF.MaCa --- preliminary version of \today\ --- Padova
%\hss
} 

\footline{
\hbox to0pt{\small University of Padova, Italy\hss}
\hss\rm\folio\hss
\hbox to0pt{\hss\small %M.Cailotto\copyright 2004---%
\smalltt maurizio@math.unipd.it}}

%%%%%%%%%%%%%%%%%%%%%%%%%%%%%%%%%%%%%

\title{Deligne Localized Functors.}

\author{Maurizio Cailotto}

%%%%%%%%%%%%%%%%%%%%%%%%%%%%%%%%%%%% 
%      abstract 
%%%%%%%%%%%%%%%%%%%%%%%%%%%%%%%%%%%%%%%%

\abstract{ 
In this paper we present the notion of 
``Deligne localized functors'', an avatar of the derived functors, 
whose definition is inspired by Deligne in [SGA 4,XVII]. 
Their definition involves the notions of 
Ind and Pro categories, they always exist and are characterized in 
terms of universal properties. 
The classical localized functor, 
in the sense of Grothendieck and Verdier, 
exists if suitable conditions are verified for the Deligne localized 
functors. 
We apply these notions to triangulated and derived categories. 
}

\footnote{}{2000 AMS Subject Classification: 18E25,30,35.}

%%%%%%%%%%%%%%%%%%%%%%%%%%%%%%%%%%%% 
\intro 
%%%%%%%%%%%%%%%%%%%%%%%%%%%%%%%%%%%%%%%%%%%

This paper arises from a study of  duality theorems in 
algebraic geometry (Grothendieck, Hartshorne, Deligne,...) 
and attempts to understand the abstract structure of 
that theory, in order to extend it to other contexts. 
To that end, a preliminary study of an abstract version of 
the notion of derived functors turns out to be indispensable.  Here  we present 
that preliminary step, with a constant attention to its generality 
and possible applications to wider contexts. 

We begin with an overview of the contents of this article. 
Paragraph $\S 0$ contains a review of notation and results 
on well known arguments: $\Ind$ and $\Pro$ categories, 
localization of categories. 
It is inserted only for ease of reference 
(which would otherwise be spread out in various papers) and to insert some 
specific points for which there seems to be no adequate reference. 
Essentially no proofs are reported, and the sources of the material 
are [SGA4], [SGA1], [BBD], [KS], \dots  
The reader is advised to skip this paragraph, and to refer to it, 
if necessary, only in reading the principal matters at hand. 

Paragraphs $\S 1$ and $\S 2$ concern the notion of derived
functor, especially the Deligne definition of derived functor as 
sketched in [SGA4,{\it XVII},\S 1]. 
The formalism and the ideas underlying this theory constitute our
guide in the sequel. 
In fact we prefer the notion of localized functors, 
because it simplifies the terminology, without loss of any 
important aspect of the theory. 
(A similar point of view is taken independently in [KS2].) 
We hope that the exposition will be useful for other
applications, and also for a better understanding of the theoretical
status of (the usual) derived functors. 

In the paragraph $\S 3$ we apply these notions to the case of 
triangulated categories. The principal problem here is that 
the categories of ind- and pro-objects of a triangulated category 
are no longer triangulated. Thus special care is need to manipulate 
the notion of distinguished triangles. 

The applications to derived categories 
are given in the last paragraph. 

 The origin and the development of this work was fostered
through many discussions with Luisa Fiorot and Francesco Baldassarri. 
I would to thank both of them for suggestions and comments on 
preliminary versions of this paper.  

The author was partially supported by the grant PGR ``CPDG021784'' 
(University of Padova, Italy) 
during the preparetion of this work. 

%%%%%%%%%%%%%%%%%%%%%%%%%%%%%%%%%%%% 
\vskip 15pt  
{\sectionfont Contents}\vskip 10pt
{\proc Introduction.\endgraf }
{\proc 0.\ Notation and Preliminaries. \endgraf }
{\proc 1.\ Deligne localized functors. \endgraf }
{\proc 2.\ Grothendieck-Verdier localized functors. \endgraf }
{\proc 3.\ The case of Triangulated Categories.\endgraf }
{\proc 4.\ The case of Derived Categories.\endgraf }
%\medskip
%%%%%%%%%%%%%%%%%%%%%%%%%%%%%%%%%%%%%%%

%%%%%%%%%%%%%%%%%%%%%%%%%%%%%%%%%%%% 
\n=-1
\section{Notation and Preliminaries. }
%\rif{SECNOTATIONS}
%%%%%%%%%%%%%%%%%%%%%%%%%%%%%%%%%%%%%%

%%%%%%%%%%%%%%%%%%%%%%%%%%%%%%%%%%%%%% 
{\numfont\NN. \bf $\Pro$ and $\Ind$ categories. }
\rif{PREPROINDCAT}

\smallskip 

{\numfont\NNN. \proc Functors Categories.} 
\rif{PREFUNCAT}%
For a category $\c C$ we put 
$\c C^\cov\bydef\Funct(\c C,\Set)$ and 
$\c C^\con\bydef\Funct(\c C^\opp,\Set)$ 
the categories of covariant and contravariant functors to 
the category of sets. 
We call $h^\cov:\c C^\opp\ra\c C^\cov$ and 
$h^\con:\c C\ra\c C^\con$ the canonical functors sending 
an object to its representable 
(covariant and contravariant) functors. 

For any $F\in\c C^\cov$ we may define the category $\c C/F$ 
(object of $\c C$ endowed with a morphism $h^\cov(X)\ra F$, 
and morphisms compatible with these data); the canonical 
morphism  
$\limin_{X\in\c C/F}h^\cov(X)\lra F$ 
(in $\c C^\cov$) is an isomorphism. 
In particular we may describe the morphisms in $\c C^\cov$ 
between two functors as 
$$ 
\eqalign{
\Hom_{\c C^\cov}(F,G)& = 
\Hom_{\c C^\cov}(\limin_{X\in\c C/F}h^\cov(X),\limin_{Y\in\c C/G}h^\cov(Y))
\cr &= 
\limpr_{X\in\c C/F}\Hom_{\c C^\cov}(h^\cov(X),\limin_{Y\in\c C/G}h^\cov(Y))
\cr &= 
\limpr_{X\in\c C/F}\limin_{Y\in\c C/G}h^\cov(Y)(X)\cr &= 
\limpr_{X\in\c C/F}\limin_{Y\in\c C/G}\Hom_{\c C}(Y,X) 
}$$ 
which is the general formulation of the Yoneda lemma 
(in the last equality the standard Yoneda lemma is used). 

Dually we have that $F\in\c C^\con$ is isomorphic to 
$\limin_{X\in\c C/F}h^\con(X)$ 
and 
$$
\Hom_{\c C^\con}(F,G)= 
\limpr_{X\in\c C/F}\limin_{Y\in\c C/G}\Hom_{\c C}(X,Y)
$$

\smallskip 

{\numfont\NNN. \proc (Pseudo-)Filtrant Categories.} 
\rif{PREFILTCAT}%
A category $\c I$ is pseudo-filtrant if the following conditions hold: 
\item{$(PF1)$}
any diagram $j\la i\ra j'$ 
% $\diagram{i &\R\w &j \cr \D\h \cr j'\cr}$ 
can be completed with 
$j\ra k\la j'$
to form a commutative square; 
% $\diagram{i &\R\w &j \cr \D\h&&\D\h \cr j'&\R\w &k \cr}$; 
\item{$(PF2)$}
any diagram $i\RR\w j$ can be completed to 
$i\RR\w j\R\w k$ 
commutative 
(any two parallel morphisms can be equalized). 
\endgraf\noindent
A non-empty, pseudo-filtrant category is filtrant if it is connected 
(i.e. any two objects can be connected by a sequence of morphisms, 
independently of the directions). 
Note that 
\item{$(i)$} 
under the condition $(PF1)$, a category is connected iff 
the following condition $(C')$ holds: 
for any $i,j\in\ob\c I$ there exists $k\in\ob\c I$ and a diagram 
$i\ra k\la j$; 
\item{$(ii)$} 
$(PF2)$ and $(C')$ imply $(PF1)$; 
\item{$(iii)$} 
in particular, $\c I$ is filtrant iff it is non empty, 
$(PF2)$ and $(C')$ hold. 
\endgraf\noindent 

In particular, a category with amalgamed sums and cokernels 
is pseudo-filtrant, 
and a category with finite sums and cokernels 
is filtrant. 

A filtrant category $\c I$ is essentially small if 
it admits a small full subcategory $\c I'$ which is cofinal, 
i.e. such that for any functor $F:\c I\ra\c C$ the inclusion 
$i:\c I'\ra\c I$ induces an isomorphism 
$\limin_{\c I'}F\circ i\la\limin_{\c I}F$ 
in $\c C^\con$; 
or equivalently such that for any $G:\c I^\opp\ra\c C$ the inclusion 
$i:\c I'\ra\c I$ induces an isomorphism 
$\limpr_{\c I'}G\circ i\la\limpr_{\c I}G$ 
in $\c C^\cov$. 

Observe that for any essentially small filtrant category $\c I$ 
there exists a small filtrant ordered set $E$ with a cofinal 
functor ${(E,\leq)}\ra\c I$ (Deligne [RD, App. n$^o$ 1]). 

\smallskip 

{\numfont\NNN. \proc Reverse Inductive Limits.} 
\rif{PREREVINDLIM}%
Let $F:\c I\ra\c C$ with $\c I$ a small filtrant category; 
we define the functor 
$\uplimin_{\c I}F:\c C^\opp\ra\c Set$ (i.e. in $\c C^\con$) 
by 
$\uplimin_{\c I}F=\limin_{\c I}h^\con(F)=\limin_{i\in\c I}h^\con(Fi)$ 
(inductive limit in the category $\c C^\con$), 
i.e. 
$(\uplimin_{\c I}F)(X)=
\limin_{\c I}\Hom_{\c C}(X,F)=\limin_{\c I}h^\cov(X)F$. 

It is representable if there exists $L\in\ob\c C$ such that 
$h^\con(L)\isom\uplimin_{\c I}F$, i.e. for any $W\in\ob\c C$ 
we have 
$\Hom_{\c C}(W,L)\isom\limin_{\c I}\Hom_{\c C}(W,F)$, 
bijection realized by the universal property of (the class of) 
a morphism $f:L\ra F(i_0)$, that is: 
for any $u:W\ra F(i)$ there exists a unique $\phi:W\ra L$ 
such that the classes of $f\phi$ and $u$ coincide, 
i.e. such that there exists $i_0\ra^{s_0} k\la^s i$ 
with $F(s_0)f\phi=F(s)u$. 

The representative $L$ is characterized by the following properties: 
\item{$(1)$} 
for any $i\in\ob\c I$ there exists $\iota_i:F(i)\ra L$; 
\item{$(2)$} 
there exists $i_0\in\ob\c I$ and a morphism $f:L\ra F(i_0)$; 
\endgraf\noindent such that 
\item{$(a)$} 
for any $i$ there exists $i_0\ra^{s_0} k\la^s i$ 
such that $F(s_0)f\iota_i=F(s)$; 
\item{$(b)$} 
$\iota_{i_0}f=\id_L$; 
\item{$(c)$} 
for any $s:i\ra j$ we have $\iota_i=\iota_jF(s)$. 
\endgraf\noindent 
In fact the bijections 
$\Hom_{\c C}(W,L)\ra\limin_{\c I}\Hom_{\c C}(W,F)$ 
are realized by sending $\phi$ to the class of $f\phi$ 
with inverse sending the class of $f_i$ to $\iota_if_i$. 

All functors $T:\c C\ra\c C'$ preserve the representative $L$ 
of $\uplimin_{\c I}F$, i.e. $T(L)$ is always a representative of 
$\uplimin_{\c I}T\circ F$. 

If $\uplimin_{\c I}F$ is represented by $L$, then also 
the functor $\limin_{\c I}F$ is represented by $L$, 
using the bijection 
$\Hom_{\c C}(L,W)\ra\limpr_{\c I}\Hom_{\c C}(F,W)$ 
sending $\phi$ to the sequence $(\phi\iota_i)$, 
and inverse sending $(f_i)$ to $f_{i_0}f$. 
Therefore, if the category admits inductive limits, 
we have necessarily $L\isom\limin_{\c I}F$, with universal data 
given by $(1)$ and $(c)$. 
The functor $\uplimin_{\c I}F$ is representable if and only if 
the canonical morphism $c:\uplimin_{\c I}F\lra\limin_{\c I}F$ 
(in $\Ind(\c C)$, see below) is an isomorphism, 
i.e. if and only if there exists an inverse 
$f:\limin_{\c I}F\lra\uplimin_{\c I}F$ 
(corresponding to $(2)$) with $cf=\id$ 
(corresponding to $(b)$) and $fc=\id$ 
(corresponding to $(a)$). 
In that case any functor $T:\c C\ra\c C'$ commutes 
with the inductive limit of the system $F$. 

Note that 
$\Hom_{\c C^\con}(\uplimin_{\c I}F,H)\isom\limpr_{\c I}H\circ F$ 
and that for $F:\c I\ra\c C$ and 
$G:\c J\ra\c C$ in $\c C^\con$ we have 
$$ 
\Hom_{\c C^\con}(\uplimin_{\c I}F,\uplimin_{\c J}G)\isom 
\limpr_{\c I}\limin_{\c J}\Hom_{\c C}(F,G) = 
\limpr_{i\in\c I}\limin_{j\in\c J}\Hom_{\c C}(F_i,G_j)\ . 
$$ 

\smallskip 

{\numfont\NNN. \proc  Ind-Objects.} 
\rif{PREINDOBJ}%
We define the category of Ind-object of $\c C$ equivalently as either 
\item{$(i)$} 
the full subcategory $\Ind\-\c C$ of $\c C^\con$ whose objects are 
the functors isomorphic to filtrant inductive limits of 
representable functors; or
\item{$(ii)$} 
the category $\Ind(\c C)$ whose objects are the filtrant inductive
systems, i.e. the functors $F: \c I\ra\c C$ from a small filtrant
category, and morphisms defined by 
$\Hom_{\Ind(\c C)}(F,G)=
\limpr_{i\in\c I}\limin_{j\in\c J}\Hom_{\c C}(Fi,Gj)$. 
\endgraf\noindent 
The equivalence $\Ind(\c C)\ra\Ind\-\c C\subseteq\c C^\con$ 
between these two categories is defined by 
sending a functor $F:\c I\ra\c C$ to $\uplimin_{\c I}F$. 

Then we have that the functor $h^\con:\c C\ra\c C^\con$ 
extends to a left exact fully faithful functor 
$h^\con:\Ind(\c C)\ra\c C^\con$ 
and restricts to an exact fully faithful functor 
$h^\con:\c C\ra\Ind(\c C)$
making commutative the diagram 
$$ 
\diagram{
\c C &\R\W^i &\Ind(\c C) \cr 
\ln{h^\con}\D\H & &\D\H\rn{h^\con} \cr 
\Ind\-\c C &\R\W_i &\c C^\con \cr
}$$ 
where the first morphism of any edge 
is an exact fully faithful functor. 

We remark that in general the canonical functors 
$\c Funct(\c I,\c C)\ra\Ind(\c C)$ 
are neither full nor faithful. 

Suppose that $\c C$ admits  filtrant inductive limits; 
then the canonical bijection 
$$\Hom_{\c C}(\limin_{\c I}F,W)\isom\Hom_{\Ind(\c C)}(F,W)$$
shows that $\limin_{\c I}$ is the left adjoint of the canonical
inclusion $\c C\ra\Ind(\c C)$. 
Moreover the following conditions are equivalent: 
\item{$(a)$} 
the canonical functor $\c C\ra\Ind(\c C)$ commutes with  
filtrant inductive limits; 
\item{$(b)$} 
for any $X\in\ob\c C$ the functor $h^\cov(X)\in\c C^\cov$ 
commutes with  filtrant inductive limits; 
\item{$(c)$} 
the functor $\limin:\Ind(\c C)\ra\c C$ is fully faithful 
(and so is an equivalence of categories). 
\endgraf 

%% In particular we have $\Ind(\Ind(\c C))\not\isom\Ind(\c C)$. 

\smallskip 

{\numfont\NNN. \proc Extension of functors.}
\rif{PREEXTFUNIND}%
Let $F:\c C\ra\Ind\c D$ be a functor. 
Then we may extend $F$ to a functor 
$\ov F:\Ind\c C\ra\Ind\c D$
uniquely defined by the condition of commutation 
with $\uplimin$, that is 
$\ov F(\uplimin_{\c I}X_{i})\bydef\uplimin_{\c I}F(X_{i})$. 
This defines a functor 
$\Funct(\c C,\Ind\c D)\lra\Funct(\Ind\c C,\Ind\c D)$ 
which is fully faithfull. 
The image of a morphism $\phi:F\ra G$ is denoted 
$\ov\phi:\ov F\ra\ov G$. 

\smallskip 

{\numfont\NNN. \proc Double $\Ind$ categories.}
\rif{PREDOUINDCAT}%
The category $\Ind(\c C)$ admits filtrant inductive limits, 
so that we have a functor 
$\uplimin:\Ind(\Ind(\c C))\lra\Ind(\c C)$ 
which is an exact left adjoint of the canonical inclusion. 
In general it is not fully faithful. 

\smallskip 

{\numfont\NNN. \proc Strict Ind-objects.}
\rif{PRESTRINDOBJ}%
An ind-object $\phi:\c I\ra\c C$ is strict if 
$\c I$ is (the category associated to) a small ordered set, 
and one of the following equivalent condition holds: 
\item{$(i)$} 
the canonical morphisms $\phi(i)\ra\uplimin_{\c I}\phi$ 
are monomorphisms in $\c C^\con$; 
\item{$(ii)$} 
for any $i\leq j$ the transition morphism $\phi(i)\ra\phi(j)$ 
is a monomorphism in $\c C$. 
\endgraf\noindent  
It is essentially strict if it is isomorphic in $\Ind(\c C)$ 
to a strict one. 

\smallskip 

{\numfont\NNN. \proc   Constant Ind-objects.} 
\rif{PRECONSINDOBJ}%
An ind-object is said to be constant if it is in the image of the
canonical functor $\c C\ra\Ind(\c C)$, 
and essentially constant if it is in the essential image, i.e. 
if is isomorphic in $\Ind(\c C)$ to a constant one. 

\smallskip 

{\numfont\NNN. \proc  Ind-Representability.} 
\rif{PREINDREPR}%
A functor $F\in\c C^\con$ is ind-representable if it is in the 
essential image of the inclusion $\Ind\-\c C\ra\c C^\con$, 
i.e. if it is isomorphic to an inductive limit in 
$\c C^\con$ of representable functors. 

An ind-representable functor $F$ is left exact, 
i.e. the canonical morphism 
$F(\limin_{\c I}\phi)\ra\limpr_{\c I}F\phi$ 
is an isomorphism for all finite projective systems 
$\phi:\c I\ra\c C$. 

\smallskip 

{\numfont\NNN. \proc  Criterion of ind-representability.} 
\rif{PREINDREPRCRI}%
The following conditions are equivalent: 
\item{$(a)$} 
$F$ is ind-representable; 
\item{$(b)$} 
the category $\c C/F$ is essentially small and filtrant; 
\item{$(b')$} 
if the category $\c C$ is equivalent to a small category: 
$\c C/F$ is filtrant; 
\item{$(c)$} 
if in $\c C$ all finite inductive limits are representable: 
$F$ is a left exact functor and $\c C/F$ is essentially small; 
\item{$(c')$} 
if the category $\c C$ is equivalent to a small category 
and in $\c C$ all finite inductive limits are representable: 
$F$ is a left exact functor. 
\endgraf\noindent 
Note that if $\c C$ has finite inductive limits, then 
$F$ left exact implies that $\c C/F$ also has finite inductive limits,
so that, in particular, it is filtrant. 

For $F\in\c C^\con$, let $Sub(F)$ be the full subcategory of 
$\c C/F$ given by the injective morphisms (i.e. the representable
sub-functors of $F$). 
Then the following are equivalent: 
\item{$(i)$} 
$F$ is strictly ind-representable (i.e. ind-representable by a strict 
ind-object); 
\item{$(ii)$} 
the category $Sub(F)$ is filtrant, essentially small and 
cofinal in $\c C/F$. 
\endgraf

\smallskip 

{\numfont\NNN. \proc  Ind-Adjoints.} 
\rif{PREINDADJ}%
Consider $F:\c C\ra\c C'$ a functor, 
and $F^\con:\c C^{\prime\con}\ra\c C^\con$ the canonical 
inverse image; 
we say that $F$ admits an ind-adjoint if one of the following 
equivalent conditions are satisfied: 
\item{$(a)$} 
$F^\con$ sends $\Ind\-\c C'$ in $\Ind\-\c C$; 
\item{$(a')$} 
$F^\con$ sends $\c C'$ in $\Ind\-\c C$; 
\item{$(b)$} 
for any $Z'\in\ob\Ind(\c C')$ the functor in $\c C^\con$ 
sending $X$ to $\Hom_{\Ind(\c C')}(FX,Z')=\limin\Hom_{\c C}(FX,Z')$  
is ind-representable; 
\item{$(b')$} 
for any $X'\in\ob\c C'$ the functor in $\c C^\con$ 
sending $X$ to $\Hom_{\c C}(FX,X')$ is ind-representable; 
\item{$(c)$} 
there exists a functor $G:\Ind\-\c C'\ra\Ind\-\c C$ such that 
we have a bifunctorial isomorphism 
$\Hom_{\Ind\-\c C}(X,GZ')\isom\Hom_{\Ind\-\c C'}(FX,Z')$ 
for any $X\in\ob\c C$ and $Z'\in\ob\Ind\-\c C'$; 
\item{$(c')$} 
there exists a functor $G_0:\c C'\ra\Ind\-\c C$ such that 
we have a bifunctorial isomorphism 
$\Hom_{\Ind\-\c C}(X,G_0X')\isom\Hom_{\c C'}(FX,X')$ 
for any $X\in\ob\c C$ and $X'\in\ob\c C'$; 
\item{$(d)$} 
the functor $\Ind(F):\Ind(\c C)\ra\Ind(\c C')$ 
admits a right adjoint; 
\item{$(e)$} 
if $\c C$ is equivalent to a small category: 
$F$ is right exact. 
\endgraf

Note that if $F$ admits a right adjoint $G':\c C'\ra\c C$, 
then it admits an ind-adjoint which is canonically 
isomorphic to $\Ind(G')$. 

\smallskip 

{\numfont\NNN. 
\proc  Presentation of morphisms of Ind-objects.}\rif{PREMORIND} 
Let $\c F\ell (\c C)$ be 
the category of morphisms of the category $\c C$ 
(morphisms of $\c F\ell (\c C)$ are the commutative squares). 
Then we have a canonical functor 
$$ 
\Ind(\c F\ell (\c C))\R\W\c F\ell (\Ind(\c C)) 
$$ 
which is fully faithful and admits a quasi inverse 
right adjoint (``parallelization of morphisms'') 
given by the following construction (Artin-Mazur). 
Let $f:X\ra Y$ be a morphism with $X$ and $Y$ objects of $\Ind(\c C)$, 
inductive systems indexed by $i\in\ I$ and $j\in J$ respectively. 
We say that a morphism $\phi:X_i\ra Y_j$ is a component of $f$ 
if $\phi=f_i\in\limin_j\Hom_{\c C}(X_i,Y_j)$, i.e. if the 
following natural diagram 
$$ 
\diagram{
X_i  &\R\w^\phi &Y_j  \cr 
\D\h &          &\D\h \cr 
X    &\R\w_f    &Y    \cr
}$$
commutes in $\Ind(\c C)$. 

Let $\Phi_f$ be the category whose objects are the components of $f$, 
and morphisms from $\phi:X_i\ra Y_j$ to $\phi':X_{i'}\ra Y_{j'}$ 
are the data of two morphisms 
$m:i\ra i'$ (in $I$) and $n:j\ra j'$ (in $J$) 
such that $Y_n\phi=\phi'X_m$ (morphisms of $\c C$). 
Then $\Phi_f$ is a small filtering category and the natural 
projection functors 
$\Phi_f\ra I$ and $\Phi_f\ra J$ are cofinal functors. 
Therefore we may associated to $f$ the object of 
$\Ind(\c F\ell (\c C))$ given by the system 
$\phi:X_\phi=X_i\ra Y_j=Y_\phi$ 
indexed by $\phi$ in the category $\Phi_f$. 
We may prove the adjoint property and the functoriality of the
construction. 
More precisely, the construction gives the Ind-adjoint 
of the canonical functor 
$\c F\ell (\c C)\R\W\c F\ell (\Ind(\c C))$.  

In the same way (see [AM]) we can prove the 
{\proc Uniform Approximation Lemma}: 
let $\Delta$ be a finite type of diagram with commutativity conditions 
without loops (i.e. a finite category without loops), 
and denote by $\Delta(\c C)$ 
the category of diagrams of type $\Delta$ in $\c C$, 
that is, the category of functors from $\Delta$ to $\c C$; 
then the natural functor 
$$ 
\Ind(\Delta(\c C))\R\W \Delta(\Ind(\c C))
$$ 
admits a right adjoint which is a quasi-inverse. 

% We remark that the canonical isomorphisms attached to the
% adjunction situation involve only identity morphisms of the
% objects in the inductive systems. 

%%%%%%%%%%%%%%%%%%%%%%%%%%%%%%%%%%%%%%%

\smallskip 

{\numfont\NNN. \proc }
The previous construction can be dualized in order to define 
the category of pro-objects as 
$\Pro(\c C)\bydef(\Ind(\c C^\opp))^\opp$. 
In the following we make explicit all the definitions and results. 

\smallskip 

{\numfont\NNN. \proc  Reverse Projective Limits.}
\rif{PREREVPROLIM}%
Let $F:\c I^\opp\ra\c C$ with $\c I$ a small filtrant category; 
we define the functor 
$\uplimpr_{\c I}F:\c C\ra\c Set$ (i.e. in $\c C^\cov$) 
by 
$\uplimpr_{\c I}F=\limin_{\c I}h^\cov(F)=\limin_{i\in\c I}h^\cov(Fi)$ 
(inductive limit in the category $\c C^\cov$), 
i.e. 
$(\uplimpr_{\c I}F)(X)=
\limin_{\c I}\Hom_{\c C}(F,X)=\limin_{\c I}h^\con(X)F$. 

This $\uplimpr$ is representable if there exists $M\in\ob\c C$ such that 
$h^\cov(M)\isom\uplimpr_{\c I}F$, i.e. for any $W\in\ob\c C$ 
we have 
$\Hom_{\c C}(M,W)\isom\limin_{\c I}\Hom_{\c C}(F,W)$, 
with the bijection being realized by the universal property of (the class of) 
a morphism $f:F(i_0)\ra M$: 
for any $u:F(i)\ra W$ there exists a unique $\phi:M\ra W$ 
such that the classes of $\phi f$ and $u$ coincide, 
i.e. such that there exists $i_0\ra^{s_0} k\la^s i$ 
with $\phi fF(s_0) =uF(s)$. 

The representative $M$ is characterized by the following properties: 
\item{$(1)$} 
for any $i\in\ob\c I$ there exists $\pi_i:M\ra F(i)$; 
\item{$(2)$} 
there exists $i_0\in\ob\c I$ and a morphism $f:F(i_0)\ra M$; 
\endgraf\noindent
such that 
\item{$(a)$} 
for any $i$ there exists $i_0\ra^{s_0} k\la^s i$ 
such that $\pi_i fF(s_0)=F(s)$; 
\item{$(b)$} 
$f\pi_{i_0}=\id_M$; 
\item{$(c)$} 
for any $s:i\ra j$ we have $\pi_i=F(s)\pi_j$. 
\endgraf\noindent 
In fact the bijections 
$\Hom_{\c C}(M,W)\ra\limin_{\c I}\Hom_{\c C}(F,W)$ 
are realized by sending $\phi$ to the class of $\phi f$ 
with inverse sending the class of $f_i$ to $f_i\pi_i$. 

All functors $T:\c C\ra\c C'$ preserve the representative $M$ 
of $\uplimpr_{\c I}F$, i.e. $T(L)$ is always a representative of 
$\uplimpr_{\c I}T\circ F$. 

If $\uplimpr_{\c I}F$ is represented by $M$, then also 
the functor $\limpr_{\c I}F$ is represented by $M$, 
using the bijection 
$\Hom_{\c C}(W,M)\ra\limpr_{\c I}\Hom_{\c C}(W,F)$ 
which sends $\phi$ to the sequence $(\pi_i\phi)$, 
with inverse sending $(f_i)$ to $ff_{i_0}$. 
Therefore, if the category admits projective limits, 
we have necessarily $L\isom\limpr_{\c I}F$, with universal data 
given by $(1)$ and $(c)$. 
The functor $\uplimpr_{\c I}F$ is representable if and only if 
the canonical morphism $c:\limpr_{\c I}F\lra\uplimpr_{\c I}F$ 
(in $\Pro(\c C)$, see below) is an isomorphism, 
i.e. if and only if there exists an inverse 
$f:\uplimpr_{\c I}F\lra\limpr_{\c I}$ 
(corresponding to $(2)$) with $fc=\id$ 
(corresponding to $(b)$) and $cf=\id$ 
(corresponding to $(a)$). 

Note that 
$\Hom_{\c C^\cov}(\uplimpr_{\c I}F,H)\isom\limpr_{\c I}H\circ F$ 
and that for $F:\c I^\opp\ra\c C$ and 
$G:\c J^\opp\ra\c C$ we have 
$$ 
\Hom_{\c C^\cov}(\uplimpr_{\c I}F,\uplimpr_{\c J}G)\isom 
\limpr_{\c I}\limin_{\c J}\Hom_{\c C}(G,F)\ . 
$$ 

\smallskip 

{\numfont\NNN. \proc  Pro-Objects.}
\rif{PREPROOBJ}%
We define the category of Pro-object of $\c C$ (anti)equivalently as: 
\item{$(i)$} 
the full subcategory $\Pro\-\c C$ of $\c C^\cov$ whose objects are 
the functors isomorphic to filtrant inductive limits of 
representable functors; 
\item{$(ii)$} 
the category $\Pro(\c C)$ whose objects are the filtrant projective
systems, i.e. the functors $F: \c I^\opp\ra\c C$ from a small filtrant
category, and morphisms defined by 
$\Hom_{\Pro(\c C)}(F,G)=
\limpr_{j\in\c J}\limin_{i\in\c I}\Hom_{\c C}(Fi,Gj)$. 
\endgraf\noindent 
The (anti)equivalence $\Pro(\c C)^\opp\ra\Pro\-\c C\subseteq\c C^\cov$ 
between the two categories is defined by 
sending a functor $F:\c I^\opp\ra\c C$ to $\uplimpr_{\c I}F$. 

Then we have that the functor $h^\cov:\c C^\opp\ra\c C^\cov$ 
extends to a left exact fully faithful functor 
$h^\cov:\Pro(\c C)^\opp\ra\c C^\cov$ 
which makes  the following diagram commutative
$$ 
\diagram{
\c C^\opp &\R\W^{i^\opp} &\Pro(\c C)^\opp \cr 
\ln{h^\cov}\D\H & &\D\H\rn{h^\cov} \cr 
\Pro\-\c C &\R\W_i &\c C^\cov \cr
}$$ 
where the first morphism of any edge is an exact fully faithful 
functor. 

Observe that in general the canonical functors 
$\c Funct(\c I^\opp,\c C)\ra\Pro(\c C)$ 
are neither full nor faithful. 

Suppose that $\c C$ admits  filtrant projective limits; 
then the canonical bijection 
$\Hom_{\c C}(W,\limpr_{\c I}F)\isom\Hom_{\Pro(\c C)}(W,F)$
shows that $\limpr_{\c I}$ is the right adjoint of the canonical
inclusion $\c C\ra\Pro(\c C)$. 
Moreover the following conditions are equivalent: 
\item{$(a)$} 
the canonical functor $\c C\ra\Pro(\c C)$ commutes with  
filtrant inductive limits; 
\item{$(b)$} 
for any $X\in\ob\c C$ the functor $h^\con(X)\in\c C^\con$ 
commutes with filtrant inductive limits; 
\item{$(c)$} 
the functor $\limpr:\Pro(\c C)\ra\c C$ is fully faithful 
(and so an equivalence of categories). 
\endgraf

%% In particular we have $\Pro(\Pro(\c C))\not\isom\Pro(\c C)$. 

\smallskip 

{\numfont\NNN. \proc Extension of functors.}
\rif{PREEXTFUNPRO}%
Let $F:\c C\ra\Pro\c D$ be a functor. 
Then we may extend $F$ to a functor 
$\ov F:\Pro\c C\ra\Pro\c D$
uniquely defined by the condition of commutation 
with $\uplimpr$, that is 
$\ov F(\uplimpr_{\c I}X_{i})\bydef\uplimpr_{\c I}F(X_{i})$. 
This defines a functor 
$\Funct(\c C,\Pro\c D)\lra\Funct(\Pro\c C,\Pro\c D)$ 
which is fully faithfull. 
The image of a morphism $\phi:F\ra G$ is denoted 
$\ov\phi:\ov F\ra\ov G$.

\smallskip 

{\numfont\NNN. \proc Double $\Pro$ categories.}
\rif{PREDOUPROCAT}%
The category $\Pro(\c C)$ admits  filtrant inductive limits, 
so that we have a functor 
$\uplimpr:\Pro(\Pro(\c C))\lra\Pro(\c C)$ 
which is an exact right adjoint of the canonical inclusion. 
In general it is not fully faithful. 

\smallskip 

{\numfont\NNN. \proc  Strict Pro-objects.}
\rif{PRESTRPROOBJ}%
A pro-object $\phi:\c I^\opp\ra\c C$ is strict if 
$\c I$ is (the category associated to) a small ordered set, 
and one of the following equivalent condition holds: 
\item{$(i)$} 
the canonical morphisms $\uplimpr_{\c I}\phi\ra\phi(i)$ 
are epimorphisms in $\c C^\cov$; 
\item{$(ii)$} 
for any $i\leq j$ the transition morphism $\phi(j)\ra\phi(i)$ 
is an epimorphism in $\c C$. 
\endgraf\noindent  
The pro-object $\phi$ is essentially strict if it is isomorphic in $\Pro(\c C)$ 
to a strict pro-object. 

\smallskip 

{\numfont\NNN. \proc   Constant Pro-objects.} 
\rif{PRECONPROOBJ}%
A pro-object is said to be constant if it is in the image of the
canonical functor $\c C\ra\Pro(\c C)$, 
and essentially constant if it is in the essential image, i.e., 
if is isomorphic in $\Pro(\c C)$ to a constant pro-object. 

\smallskip 

{\numfont\NNN. \proc   Pro-Representability.} 
\rif{PREPROREPR}%
A functor $F\in\c C^\cov$ is pro-representable if it is in the 
essential image of the inclusion $\Pro\-\c C\ra\c C^\cov$, 
i.e., if it is isomorphic to an inductive limit in 
$\c C^\cov$ of representable functors. 

A pro-representable functor $F$ is left exact, 
i.e. the canonical morphism 
$F(\limpr_{\c I}\phi)\ra\limpr_{\c I}F\phi$ 
is an isomorphism for all finite projective system 
$\phi:\c I\ra\c C$. 

\smallskip 

{\numfont\NNN. \proc  Criterion of pro-representability.} 
\rif{PREPROREPRCRI}%
The following conditions are equivalent: 
\item{$(a)$} 
$F$ is pro-representable; 
\item{$(b)$} 
the category $\c C/F$ is essentially small and filtrant; 
\item{$(b')$} 
if the category $\c C$ is equivalent to a small category: 
$\c C/F$ is filtrant; 
\item{$(c)$} 
if in $\c C$ the finite projective limits are representable: 
$F$ is a left exact functor and $\c C/F$ is essentially small; 
\item{$(c')$} 
if the category $\c C$ is equivalent to a small category 
and in $\c C$ the finite projective limits are representable: 
$F$ is a left exact functor. 
\endgraf

We remark that if $\c C$ has finite projective limits, then 
$F$ left exact implies that $\c C/F$ also has finite inductive limits,
so that, in particular, it is filtrant. 

For $F\in\c C^\cov$, let $Sub(F)$ be the full subcategory of 
$\c C/F$ given by the injective morphisms (i.e. the representable
sub-functors of $F$). 
Then the following are equivalent: 
\item{$(i)$} 
$F$ is strictly pro-representable (i.e. pro-representable by a strict 
pro-object); 
\item{$(ii)$} 
the category $Sub(F)$ is filtrant, essentially small and 
cofinal in $\c C/F$. 
\endgraf

\smallskip 

{\numfont\NNN. \proc   Pro-Adjoints.} 
\rif{PREPROADJ}%
Consider $F:\c C\ra\c C'$ a functor, 
and $F^\cov:\c C^{\prime\cov}\ra\c C^\cov$ the canonical 
inverse image; 
we say that $F$ admits a pro-adjoint if one of the following 
equivalent conditions are satisfied: 
\item{$(a)$} 
$F^\cov$ sends $\Pro\-\c C'$ in $\Pro\-\c C$; 
\item{$(a')$} 
$F^\cov$ sends $\c C'$ in $\Pro\-\c C$; 
\item{$(b)$} 
for any $Z'\in\ob\Pro(\c C')$ the functor in $\c C^\cov$ 
sending $X$ to 
$\Hom_{\Ind(\c C')}(Z',FX)=\limin\Hom_{\c C}(Z',FX)$  
is pro-representable; 
\item{$(b')$} 
for any $X'\in\ob\c C'$ the functor in $\c C^\cov$ 
sending $X$ to $\Hom_{\c C}(X',FX)$ is pro-representable; 
\item{$(c)$} 
there exists a functor $G:\Pro\-\c C'\ra\Pro\-\c C$ such that 
we have a bifunctorial isomorphism 
$\Hom_{\Pro\-\c C}(GZ',X)\isom\Hom_{\Pro\-\c C'}(Z',FX)$ 
for any $X\in\ob\c C$ and $Z'\in\ob\Ind\-\c C'$; 
\item{$(c')$} 
there exists a functor $G_0:\c C'\ra\Pro\-\c C$ such that 
we have a bifunctorial isomorphism 
$\Hom_{\Pro\-\c C}(G_0X',X)\isom\Hom_{\c C'}(X',FX)$ 
for any $X\in\ob\c C$ and $X'\in\ob\c C'$; 
\item{$(d)$} 
the functor $\Pro(F):\Pro(\c C)\ra\Pro(\c C')$ 
admits a left adjoint; 
\item{$(e)$} 
if $\c C$ is equivalent to a small category: 
$F$ is left exact. 
\endgraf

Remark that if $F$ admits a left adjoint $G':\c C'\ra\c C$, 
then it admits a pro-adjoint which is canonically 
isomorphic to $\Pro(G')$. 

\smallskip 

{\numfont\NNN. \proc   Presentation of morphisms of Pro-objects.}
\rif{PREMORPRO}%
Let again $\c F\ell (\c C)$ be 
the category of morphisms of the category $\c C$ 
(morphisms of $\c F\ell (\c C)$ are the commutative squares). 
Then we have a canonical functor 
$$ 
\Pro(\c F\ell (\c C))\R\W\c F\ell (\Pro(\c C)) 
$$ 
which is fully faithful and admits a quasi inverse 
left adjoint (``parallelization of morphisms'') 
given by the dual construction to that of \cite{PREMORIND}. 

Also the {\proc Uniform Approximation Lemma} holds: 
let $\Delta$ be a finite type of diagram with commutativity conditions 
without loops (i.e. a finite category without loops), 
and let $\Delta(\c C)$ denote 
the category of diagrams of type $\Delta$ in $\c C$, 
that is the category of functors from $\Delta$ to $\c C$; 
then the natural functor 
$$ 
\Pro(\Delta(\c C))\R\W \Delta(\Pro(\c C))
$$ 
admits a left adjoint which is a quasi-inverse.

\smallskip 

{\numfont\NNN. \proc  } \rif{PREPROINDPRO}
Occasionally we will need  categories such as $\Pro\,\Ind\c C$ 
or $\Ind\,\Pro\c C$. 
We remark only that the $\Hom_{\c C}$ as a bifunctor on 
$\c C^\opp\times\c C$ can be extended to a bifunctor 
$$ 
\Ind\,\Hom_{\c C}:(\Pro\c C)^\opp\times\Ind\c C\R\W 
\Ind\Set 
$$ 
as 
$\Ind\,\Hom_{\c C}((X_i),(Y_j))= 
\uplimin_{i\in I}\uplimin_{j\in J}\Hom_{\c C}(X_i,Y_j) 
$ 
and its composition with the inductive limit functor of $\Set$, 
i.e. 
$
\limin_{i\in I}\limin_{j\in J}\Hom_{\c C}(X_i,Y_j) 
$ 
is just the restriction of $\Hom_{\Pro\,\Ind\c C}$ 
to the (full) subcategories $\Pro\c C$ 
(first argument) and $\Ind(\c C)$ (second argument). 

Similar remarks hold for the bifunctor $\Pro\,\Hom_{\c C}$ 
with respect to the category $\Ind\,\Pro\c C$. 

In particular, we remark that for $(X_i)_{i\in I}$ in $\Pro\c C$ 
and $(Y_j)_{j\in J}$ in $\Ind\c C$ we have the equalities 
$$ 
\Hom_{\Ind\,\Pro\c C}((X_i),(Y_j))= 
\limin_{i\in I}\limin_{j\in J}\Hom_{\c C}(X_i,Y_j)= 
\limin_{j\in J}\limin_{i\in I}\Hom_{\c C}(X_i,Y_j)= 
\Hom_{\Pro\,\Ind\c C}((X_i),(Y_j)) 
$$ 
and 
$$ 
\Hom_{\Ind\,\Pro\c C}((Y_j),(X_i))= 
\limpr_{j\in J}\limpr_{i\in I}\Hom_{\c C}(Y_j,X_i)= 
\limpr_{i\in I}\limpr_{j\in J}\Hom_{\c C}(Y_j,X_i)= 
\Hom_{\Pro\,\Ind\c C}((Y_j),(X_i)) \ .
$$

\smallskip 

{\numfont\NNN. \proc  Generalized Adjunctions.} 
\rif{PREGENADJ}%
We say that 
$F:\c C\ra\Pro\-\c C'$ and $G:\c C'\ra\Ind\-\c C$ 
are generalized adjoints if 
there is a bifunctorial isomorphism 
$$\Hom_{\Pro\-\c C'}(FX,X')\isom\Hom_{\Ind\-\c C}(X,GX')$$ 
for any $X\in\ob\c C$ and $X'\in\ob\c C'$. 
In that case for any $X\in\ob\Pro\-\c C$ and $X'\in\ob\Ind\-\c C'$ 
we have  bifunctorial isomorphisms 
$$\Ind\,\Hom_{\c C'}(FX,X')\isom\Ind\,\Hom_{\c C}(X,GX')$$ 
and 
$$\Hom_{\Pro\,\Ind\-\c C'}(FX,X')\isom\Hom_{\Ind\,\Pro\-\c C}(X,GX')\ .$$ 
Moreover each of the two functors $F$ and $G$ determines the other, 
up to isomorphisms. 

We also have  the dual notions of  generalized coadjoint functors: 
$F:\c C\ra\Ind\-\c C'$ and $G:\c C'\ra\Pro\-\c C$ 
are generalized coadjoint if 
there is a bifunctorial isomorphism 
$$\Hom_{\Ind\-\c C'}(FX,X')\isom\Hom_{\Pro\-\c C}(X,GX')$$ 
for any $X\in\ob\c C$ and $X'\in\ob\c C'$. 
In that case for any $X\in\ob\Ind\-\c C$ and $X'\in\ob\Pro\-\c C'$ 
we have  bifunctorial isomorphisms 
$$\Pro\,\Hom_{\c C'}(FX,X')\isom\Pro\,\Hom_{\c C}(X,GX')$$ 
and 
$$\Hom_{\Ind\,\Pro\-\c C'}(FX,X')\isom\Hom_{\Pro\,\Ind\-\c C}(X,GX')\ .$$ 
Moreover each of the two functors $F$ and $G$ determines the other, 
up to isomorphisms. 
Note however that the notion of coadjointness is not useful 
when the categories admit (filtrant) inductive and projective limits, 
since it then reduces to the usual notion of adjunction.

\smallskip 

{\numfont\NNN. \proc  Intersection of $\Ind$ and $\Pro$.} 
\rif{PREINTINDPRO}%
The canonical square 
$$ 
\diagram{
\c C &\R\W &\Pro(\c C) \cr 
\D\H & &\D\H \cr 
\Ind(\c C) &\R\W &\Pro(\Ind(\c C)) \cr 
}$$ 
is cartesian in the following sense:  
an object of $\Pro(\Ind(\c C))$ which is in the image 
either of $\Pro(\c C)$ or of $\Ind(\c C)$ is really in $\c C$; 
i.e. $\Pro(\c C)\cap\Ind(\c C)$ in $\Pro(\Ind(\c C))$ 
is just $\c C$.

%%%%%%%%%%%%%%%%%%%%%%%%%%%%%%%%%%%%%%%%%%%%%%%%%%%%%%%%%%%%%%%%
\medskip 
{\numfont\NN. \bf Multiplicative systems and localization. }
\rif{PREMULTSYST}
\smallskip 
{\numfont\NNN. \proc right and left multiplicative systems.} 
\rif{PRELRMS}
Let $\c C$ be a category; a multiplicative system in $\c C$ 
is a family $S$ of morphism of $\c C$ such that: 
\item{$(S1)$} 
$\id_X\in S$ for any $X$ in $\c C$; 
\item{$(S2)$} 
if $f,g\in S$ then $g\circ f\in S$ if it exists. 
\endgraf\noindent  
A multiplicative system is said to be right (resp. left) 
if the following conditions are satisfied: 
\item
{$(S3)$} we may complete any diagram with  $s\in S$
$$
\diagram{&&Z\cr &&\D\h\rn{s}\cr X &\R\W_f &Y\cr} 
\qquad\hbox{in}\qquad 
\diagram{W&\R\W^g &Z\cr \ln{t}\D\h&&\D\h\rn{s}\cr X &\R\W_f &Y\cr} 
$$
commutative with $t \in S$; 
(resp. dually with the arrows reversed); 
\item
{$(S4)$} for any two morphisms $f,g:X \ra Y$ in $\c C$ 
consider the following conditions: 
\itemitem {$(i)$} there exists $s\in S$, 
                 $s:W\ra X$ such that $f\circ s=g\circ s$; 
\itemitem {$(ii)$} there exists $t\in S$, 
                $t:Y\ra Z$ such that $t\circ f=t\circ g$;
\item{} then $(i)$ implies $(ii)$ 
(resp. $(ii)$ implies $(i)$)
\endgraf\noindent   
A right and left multiplicative system is said to be bilateral, 
or simply a multiplicative system, if there is no possibility of confusion.

\smallskip 
{\numfont\NNN. \proc quasi-saturated multiplicative systems.} 
\rif{PREQSMS}
A multiplicative system $S$ is right (resp. left) 
quasi-saturated if the following 
condition holds: 
if $g\circ f\in S$ and $f\in S$ then $g\in S$ 
(resp. if $g\circ f\in S$ and $g\in S$ then $f\in S$). 
A right and left quasi-saturated multiplicative system 
is said to be quasi-saturated. 
\endgraf 

\smallskip 
{\numfont\NNN. \proc Localized categories.} 
\rif{PRELOCCAT}
Let $S$ be a multiplicative system in $\c C$; we define the 
localized category of $\c C$ with respect to $S$, 
denoted by $\c C_S$ or $\c C[S^{-1}]$, to be
a category endowed with a functor $Q:\c C \rightarrow \c C_S$ 
such that for any $s \in S$ the image $Q(s)$ is an isomorphism, 
and which is universal for this property: 
for any category $\c D$ with a functor $F:\c C \ra \c D$ 
such that $F(s)$ is isomorphism for any $s\in S$, then there exists 
a unique factorization:  
\endgraf
$$ 
\diagram{
\c C  &\R\w^Q  &\c C_S \cr 
      &\ln{F}\searrow &\D\h\rn{\exists !} \cr 
      &&\c D \ . \cr
}$$ 
More generally we may solve the same universal problem for any 
family of morphisms of a category, but the localized category has 
a rather complicated construction. 
\endgraf 
If $S$ is a right (resp. left) multiplicative system in $\c C$, then 
the category $\c C_S$ is constructed in the following way: 
the objects of $\c C_S$ are the objects of $\c C$; 
the morphisms are defined by 
$$\Hom_{\c C_S}(X,Y)=\limin_{(Y',t')\in Y/ S}\Hom_{\c C}(X,Y')$$ 
where $Y/ S$ is the full subcategory of 
$Y/\c C$ given by the objects $t':Y\ra Y'$ with $t'\in S$ 
(resp. 
$$\Hom_{\c C_S}(X,Y)=\limin_{(X',s') \in S/X}\Hom_{\c C}(X',Y)$$ 
where $S/X$ is the full subcategory of $\c C /X$ 
given by the objects 
$s':X'\ra X$, with $s'\in S$). 
The functor $Q$ is defined in the obvious way. 
\endgraf 
If the multiplicative system is bilateral, the more symmetric formula 
$$\Hom_{\c C_S}(X,Y)=
\limin_{(X',s') \in S/X\atop (Y',t')\in Y/S}
\Hom_{\c C}(X',Y')$$ 
also works. 
\endgraf 

\smallskip 
{\numfont\NNN. \proc Localization of triangulated categories.} 
\rif{PRELOCTRICAT}
Let $(\c T,T)$ be a triangulated category; 
a null system $\c N$ of $\c T$ is a family of objects 
of $\c T$ such that: 
\item{$(N1)$} $0\in\c N$,
\item{$(N2)$} $N\in\c N$ if and only if $TN\in\c N$,
\item{$(N3)$} if $X\ra Y \ra Z \ra TX$ is a distinguished triangle 
            and $X,Y\in\c N$, then $Z\in\c N$.
\endgraf 
Notice that the shift property of $\c T$ and $(N2)$ 
permit us to extend $(N3)$: 
if two vertices of a distinguished triangle are in $\c N$, 
then so is the third vertex. 
\endgraf 
In the triangulated category we can localize with respect to 
a null system; in fact the family of morphisms 
$$S(\c N)=
\l\{f:X\ra Y|\exists\hbox{ dist.tr.}X\ra^f Y\ra N\ra TX
             \hbox{ with }N\in\c N\r\}$$ 
is a quasi-saturated multiplicative system in $\c T$.  
Moreover $S(\c N)$ satisfies the following two properties: 
\item
{$(ST1)$} $s\in S(\c N)$ if and only if $T(s)\in S(\c N)$; 
\item
{$(ST2)$} if two arrows of a morphism of 
distinguished triangles are in $S(\c N)$, 
then so too is the third arrow. 
\endgraf 
{\proc Universal property:}
put $\c T/\c N=\c T_{S(\c N)}$ and let 
$Q:\c T\lra\c T/\c N$ be the canonical functor; 
then $\c T/\c N$ is canonically a triangulated category; 
$Q(N)\isom 0$ for any $N\in \c N$ and 
$\c T/\c N$ is universal with respect to this this property in the category 
of triangulated categories. 
\endgraf 
{\proc Example:}
If $\c A$ is an abelian category and 
$H:\c T\ra\c A$ is a cohomological functor, then 
the class 
$\c N_H=\{X\in\c T | H(T^nX)=0\ \forall n\in\b Z\}$ 
of $H$-acyclic objects is a null system. 
Remark that $S(\c N_H)=\{f|H(T^nf)\hbox{ iso }\forall n\in\b Z \}$.

%%%%%%%%%%%%%%%%%%%%%%%%%%%%%%%%%%%% 
\section{Deligne localized functors. }
\rif{SECDELIGNE}
%%%%%%%%%%%%%%%%%%%%%%%%%%%%%%%%%%%%%%%%%%%%%%%%%%%%%%%%%%

\proclaim{\numfont\NN. \proc Definition (Localizing functors)}. {
Let $\c C$ be a category 
and $S$ a quasi-saturated right multiplicative system in
$\c C$. Then we define the right localizing functor with respect to $S$ as 
$r'_S:\c C\lra\Ind\c C$ by 
$$r'_S(X)\bydef\uplimin_{s:X{\ra}X'}X'$$ 
where the index category is $X/S$ 
(morphisms in $S$ with source $X$), as an inductive system, 
that is 
$$r'_S(X)(Z)\bydef\limin_{s:X{\ra}X'}\Hom_{\c C}(Z,X')$$ 
as functor 
$\c C^\opp\ra\Set$. 
\endgraf 
Dually, if $S$ a quasi-saturated left multiplicative system in
$\c C$, we define the left localizing functor with respect to $S$ as 
$l'_S:\c C\lra\Pro\c C$ by 
$$l'_S(X)\bydef\uplimpr_{s:X'{\ra}X}X'$$ 
where the index category is $S/X$ 
(morphisms in $S$ with target $X$), as a projective system, 
that is 
$$l'_S(X)(Z)\bydef\limin_{s:X'{\ra}X}\Hom_{\c C}(X',Z)$$ 
as functor 
$\c C\ra\Set$. 
}\par  
\rif{DEFLOCALIZING}

{\numfont\NNN. \proc Action on morphisms.}\rif{ACTMOR}
If $f:X\ra Y$ is a morphism in $\c C$, then the morphism 
$r'_S(f):r'_S(X)\ra r'_S(Y)$ is defined in 
$$ 
\Hom_{\Ind\c C}(r'_S(X),r'_S(Y))= 
\limpr_{s:X{\ra}X'}\limin_{t:Y{\ra}Y'}\Hom_{\c C}(X',Y') 
$$ 
by the following construction: 
for any $s:X\ra X'$ we complete the diagram with $f:X\ra Y$ 
to a square 
$$ 
\diagram{
X &\R\w^s & X' \cr 
\ln{f}\D\h & &\D\h\rn{f'} \cr 
Y &\R\w_t & Y' \cr 
}$$ 
with $t\in S$; then $(t,f')$ is a representative in 
$\limin_{t:Y{\ra}Y'}\Hom_{\c C}(X',Y')$ of the 
component of $r'_S(f)$ in $s$. 
That the definition is well-posed, 
that is independent of the choices made in completing 
the square, follows from the properties of right saturated 
multiplicative systems. 
\endgraf 
>From the point of view of functors, the construction gives 
a morphism $r'_S(f):r'_S(X)\ra r'_S(Y)$ whose evaluation in $Z$ is 
$$ 
r'_S(f)(Z): 
r'_S(X)(Z)= \limin_{s:X{\ra}X'}\Hom_{\c C}(Z,X') \R\W 
\limin_{t:Y{\ra}Y'}\Hom_{\c C}(Z,Y') =r'_S(Y)(Z)
$$ 
which sends $(s,\phi)$ with $\phi:Z\ra X'$ to 
$(t,f'\phi)$, a well-defined element in the inductive limit 
on the right hand side. 
\endgraf 
Dually, if $f:X\ra Y$ is a morphism in $\c C$, then the morphism 
$l'_S(f):l'_S(X)\ra l'_S(Y)$ is defined in 
$$ 
\Hom_{\Pro\c C}(l'_S(X),l'_S(Y))= 
\limpr_{t:Y'{\ra}Y}\limin_{s:X'{\ra}X}\Hom_{\c C}(X',Y') 
$$ 
by the following construction: 
for any $t:Y'\ra Y$ we complete the diagram with $f:X\ra Y$ 
to a square 
$$ 
\diagram{
X &\L\w^s & X' \cr 
\ln{f}\D\h & &\D\h\rn{f'} \cr 
Y &\L\w_t & Y' \cr 
}$$ 
with $s\in S$; then $(s,f')$ is a representative in 
$\limin_{s:X'{\ra}X}\Hom_{\c C}(X',Y')$ of the 
component of $l'_S(f)$ in $t$. 

>From the point of view of functors, the construction gives 
a morphism $l'_S(f):l'_S(X)\ra l'_S(Y)$ whose evaluation in $Z$ is 
$$ 
l'_S(f)(Z): 
l'_S(Y)(Z)= \limin_{t:Y'{\ra}Y}\Hom_{\c C}(Y',Z)
\R\W 
\limin_{s:X'{\ra}X}\Hom_{\c C}(X',Z)  =l'_S(X)(Z)
$$ 
which sends  $(t,\psi)$ with $\psi:Y'\ra Z$ to 
$(s,\psi f')$, a well-defined element in the inductive limit.

\lemm{NNN}{
Suppose $f:X\ra Y$ in $\c C$ is a morphism in $S$, 
then $r'_{S}(f):r'_{S}(X)\ra r'_{S}(Y)$ is an isomorphism in $\Ind\c C$. 
Dually, $l'_{S}(f):l'_{S}(X)\ra l'_{S}(Y)$ is an isomorphism in $\Pro\c C$. 
}\rif{LEMINVMOR}

{\proc Proof.}
In fact we can define the inverse morphism 
$s(f):r'_{S}(Y)\ra r'_{S}(X)$ using the following construction: 
for any $t:Y\ra Y'$ in $S$, the composition with $f$ gives 
$tf:X\ra Y'$ in $S$, and this is the component of $s(f)$ in $t$. 
In terms of functor morphisms, this is the morphism 
$s(f)(Z):r'_{S}(Y)(Z)\ra r'_{S}(X)(Z)$ sending $(t,\psi)$ to $(tf,\psi)$. 
It is well defined, and the compositions with $r_{S}(f)$ are clearly 
the identities 
(one composition is easy, 
the other requires  the quasi-saturatedness 
of the multiplicative system). %need of right quasi-saturatedness 
\nobreak\hfill$\square$

\prop{NNN}{
In particular we can extend 
$r'_S:\c C\lra\Ind\c C$ to a functor 
$r_S:\c C_S\lra\Ind\c C$ and we have a commutative 
diagram of functors 
$$ 
\diagram{
\c C &\R\w^{r'_S} &\Ind\c C \cr 
\ln{Q}\D\h &\nea\rn{r_S} \cr 
\c C_S &.\cr 
}$$  
Moreover the functor $r_S$ is the ind-adjoint of the 
canonical functor $Q:\c C_S\ra\c C$. 
\endgraf 
Dually, we can extend $l'_S$ to a functor 
$l_S:\c C_S\R\W\Pro\c C$ and we have a commutative 
diagram of functors 
$$ 
\diagram{
\c C &\R\w^{l'_S} &\Pro\c C \cr 
\ln{Q}\D\h &\nea\rn{l_S} \cr 
\c C_S &.\cr 
}$$  
Moreover the functor $l_S$ is the pro-adjoint of the 
canonical functor $Q:\c C_S\ra\c C$. 
}

{\proc Proof.} 
The first claim is an obvious consequences of the lemma and the 
universal property of the localized category. 
The last claim is a consequence of the bijections 
$$ 
\Hom_{\c C_S}(QX,Y)=\limin_{t:Y{\ra}Y'}\Hom_{\c C}(X,Y')= 
\Hom_{\Ind\c C}(iX,\uplimin_{t:Y{\ra}Y'}Y')= 
\Hom_{\Ind\c C}(iX,r_SY) \ .
$$ 
The dual assertion is expressed by the bijection 
$$ 
\Hom_{\c C_S}(Y,QX)=\limpr_{s:Y'{\ra}Y}\Hom_{\c C}(Y',X)= 
\Hom_{\Pro\c C}(\uplimpr_{s:Y'{\ra}Y}Y',iX)= 
\Hom_{\Pro\c C}(l_SY,iX) \ .
$$ 
\nobreak\hfill$\square$

\defi{NNN}{
Let $i:\c C\lra\Ind\c C$ be the canonical fully faithful functor; 
then we have a natural morphism of functors 
$\delta_S:i\ra r_SQ=r'_S$, because for any $X$ 
in $\c C$ we have $\id_X\in S$. 
The morphism corresponds to the identity of $QX$ 
under the bijection of ind-adjointness between 
$Q$ and $r_S$. 
\endgraf 
Dually, let $i:\c C\lra\Pro\c C$ be the canonical fully faithful functor; 
then we have a natural morphism of functors 
$\sigma_S:l_SQ\ra i$, corresponding to the identity of $QX$ 
under the bijection of pro-adjointness between 
$Q$ and $l_S$. 
} 

For any $X$ object of $\c C$, the morphism 
$\delta_S(X):i(X)\ra r'_S(X)$ 
is represented by the identity of $X$ as a morphism between 
ind-objects, while as functors it is identified as the canonical 
morphism 
$$\delta_S(X)(Z):i(X)(Z)=\Hom_{\c C}(Z,X)\lra
\limin_{s:X{\ra}X'}\Hom_{\c C}(Z,X')=r'_S(X)(Z)$$
since the left hand side appears in the inductive limit of 
the right hand side.

\defi{NNN}{
We say that an object $X$ of $\c C$ is inert for $r_S$ or 
right inert for $S$ if $\delta_S(X)$ is an isomorphism; 
it is right localizable with respect to $S$ if $r'_S(X)$ 
is representable. 
Put $R_{S}(X)\bydef\limin_{X/S}X'$; 
then $X$ is  right localizable with respect to $S$ if and only if 
the canonical morphism $r'_S(X)\ra R_{S}(X)$ is an isomorphism 
in $\Ind\c C$. 
\endgraf 
Dually, we say that an object $X$ of $\c C$ is inert for $l_S$ or 
left inert for $S$ if $\sigma_S(X)$ is an isomorphism; 
it is left localizable with respect to $S$ if $l'_S(X)$ 
is representable. 
Put $L_{S}(X)\bydef\limin_{S/X}X'$; 
then $X$ is  left localizable with respect to $S$ if and only if 
the canonical morphism $L_S(X)\ra l'_{S}(X)$ is an isomorphism 
in $\Pro\c C$. 
} \rif{DEFLOCOBJ}\rif{DEFINERTOBJ}

The criterion \cite{PREREVINDLIM} says that 
$X$ is right localizable with respect to $S$ if and only if there 
exists $s_{0}:X\ra X_{0}$ in $S$ and a morphism 
$t_{0}:R_{S}(X)\ra X_{0}$ such that 
$\iota_{0}t_{0}=\id_{R_{S}(X)}$ 
and for any $s:X\ra X'$ there exists an object $X''$ 
in $S/X$ with morphisms 
$f':X'\ra X''$ and $f_{0}:X_{0}\ra X''$ 
in $S/X$ such that $f_{0} t_{0}\iota'=f'$ 
($\iota_{0}$and $\iota'$ indicate the canonical morphisms 
from $X_{0}$ and $X'$ to $R_{S}(X)$). 
\endgraf 
Moreover $X$ is right inert for $S$ if and only if 
$X$ is right localizable with respect to $S$ and the canonical 
morphism $X\ra R_{S}(X)$ is an isomorphism; 
that is, for any $s:X\ra X'$ in $X/S$, there exists a morphism 
$t_{X'}:X'\ra X$ such that: 
for any morphism $f:X'\ra X''$ in $X/S$ we have 
$t_{X'} f=t_{X''}$; 
$t_{X}=\id_{X}$; 
for any $s:X\ra X'$ in $X/S$ there exists $s':X\ra X''$ in $X/S$ 
and a morphism $f:X'\ra X''$ in $X/S$ such that 
$s' t_{X'}=f$.

\proclaim%
{\numfont\NNN. 
\proc Theorem (Universal property of localizing functors)}. {
For any $G:\c C_S\lra\Ind\c C$ 
the map $\beta$ to $(\beta\bullet Q)\circ\delta_{S}$ 
induces a bijection 
$$ 
\Hom_{{\Funct}(\c C_S,\Ind\c C)}(r_{S},G)\R\W
\Hom_{{\Funct}(\c C,\Ind\c C)}(i,GQ)\ ,
$$ 
that is, the pair $(r_{S},\delta_{S})$ represents the 
functor 
${\Funct}(\c C,\Ind\c C)\lra\Set$ 
which sends $G$ to $\Hom_{{\Funct}(\c C,\Ind\c C)}(i,GQ)$. 
\endgraf 
Dually, for any $G:\c C_S\lra\Pro\c C$ 
the map $\beta$ to $\sigma_{S}\circ(\beta\bullet Q)$ 
induces a bijection 
$$ 
\Hom_{{\Funct}(\c C_S,\Pro\c C)}(G,l_{S})\R\W
\Hom_{{\Funct}(\c C,\Pro\c C)}(GQ,i)\ .
$$ 
}\par\rif{LOCALIZINGUNIPRO} 

{\proc Proof.} 
The proof is a consequence of the following 
more general proposition, applied to $F=\id_{\Ind\c C}$. 
We only notice explicitly that the bijection is the composite of 
the usual one 
$$ 
\Hom_{{\Funct}(\c C_S,\Ind\c C)}(r_{S},G)\R\W
\Hom_{{\Funct}(\c C,\Ind\c C)}(r_{S}Q,GQ) 
$$ 
induced by the horizontal composition with $\id_{Q}$ 
($\beta$ to $\beta\bullet Q$) 
and the map 
$$ 
\Hom_{{\Funct}(\c C,\Ind\c C)}(r'_{S},H)\R\W
\Hom_{{\Funct}(\c C,\Ind\c C)}(i,H)\ .
$$ 
induced by the vertical composition with $\delta_{S}$ 
($\gamma$ to $\gamma\circ\delta_{S}$), which is a bijection 
if the functor $H$ sends $S$ to isomorphisms of $\Ind\c C$. 
\nobreak\hfill$\square$

\prop{NNN}{
For any functor $F:\Ind\c C\lra\Ind\c D$, 
define $r_{S}(F)\bydef Fr_{S}$. 
Then for any $G:\c C_{S}\ra\Ind\c D$ the map 
$$ 
\Hom_{{\Funct}(\c C_S,\Ind\c D)}(r_{S}(F),G)\R\W
\Hom_{{\Funct}(\c C,\Ind\c D)}(Fi,GQ) 
$$ 
induced by $\beta$ to $(\beta\bullet Q)\circ\delta_{S}(F)$, 
where $\delta_{S}(F)=F\bullet\delta_{S}$, 
is a bijection. 
\endgraf
Dually, for any functor $F:\Pro\c C\lra\Pro\c D$, 
define $l_{S}(F)\bydef Fl_{S}$. 
Then for any $G:\c C_{S}\ra\Pro\c D$ the map 
$$ 
\Hom_{{\Funct}(\c C_S,\Pro\c D)}(G,l_{S}(F))\R\W
\Hom_{{\Funct}(\c C,\Pro\c D)}(GQ,Fi) 
$$ 
induced by $\beta$ to $\sigma_{S}(F)\circ(\beta\bullet Q)$, 
where $\sigma_{S}(F)=F\bullet\sigma_{S}$, 
is a bijection. 
}\rif{GENLOCALIZINGUNIPRO}  

{\proc Proof.} 
Note that the map is the composite of 
the usual bijection  
$$ 
\Hom_{{\Funct}(\c C_S,\Ind\c D)}(r_{S}(F),G)\R\W
\Hom_{{\Funct}(\c C,\Ind\c D)}(r_{S}(F)Q,GQ) 
$$ 
induced by the horizontal composition with $\id_{Q}$ 
(sending $\beta$ to $\beta\bullet Q$),  
and the map 
$$ 
\Hom_{{\Funct}(\c C,\Ind\c D)}(r'_{S}(F),H)\R\W
\Hom_{{\Funct}(\c C,\Ind\c D)}(Fi,H)\ .
$$ 
induced by the vertical composition with $\delta_{S}(F)$ 
(sending $\gamma$ to $\gamma\circ\delta_{S}(F)$), 
which we will prove to be a bijection 
if the functor $H$ sends $S$ to isomorphisms of $\Ind\c D$ 
(as is the case for $H=GQ$ since any morphism of $S$ becames an isomorphism
after application of $Q$). 
In fact consider a morphism $\alpha:Fi\lra H$; 
by hypothesis for any 
$X$ object of $\c C$ and any $s:X\ra X'$ in $S$, we have that 
$H(s):H(X)\ra H(X')$ is an isomorphism, therefore the 
$\{H(s)^{-1}\alpha(X'):Fi(X')\ra H(X)\}$ 
is a compatible system of morphisms: 
$$ 
\diagram{
H(X') &\L\w^{\alpha(X')} &Fi(X') &\R\W^{\delta_{S}(F)(X')} &r'_S(F)(X')\cr 
\ln{H(s)}\U\h\rn{\isom} && \U\h\rn{Fi(s)} 
                          &&\ln{\isom}\U\h\rn{r'_{S}(F)(s)} \cr 
H(X) &\L\w_{\alpha(X)} &Fi(X) &\R\w_{\delta_S(F)(X)} &r'_S(F)(X) \ .\cr
}$$ 
The definition of $r'_{S}(F)(X)=\uplimin Fi(X')$ 
then gives a canonical morphism 
$r'_{S}(F)(X)\ra H(X)$ which uniquely factorizes  the given system 
through the $\{\delta_S(F)(X')\}$. 
\nobreak\hfill$\square$

\medskip 
\proclaim{\numfont\NNN. 
\proc Corollary (Localizations as adjoint functors)}. {
Let $S$ be a saturated multiplicative system in a category $\c C$; 
for any category $\c D$, consider the functor 
$$ 
Q:\Funct(\c C_S,\Ind\c D)\R\W\Funct(\c C,\Ind\c D)
$$ 
given by the composition with the canonical $Q:\c C\ra\c C_S$. 
\endgraf\noindent 
Then we can define a left adjoint 
$$r_S:\Funct(\c C,\Ind\c D)\lra\Funct(\c C_S,\Ind\c D)$$ 
of $Q$ in the following way. 
For any $F:\c C\ra\Ind\c D$ let 
$\overline{F}:\Ind\c C\ra\Ind\c D$ be its canonical extension 
to $\Ind\c C$, and define $r_S(F)\bydef\overline{F}\circ r_S$. 
Then \cite{GENLOCALIZINGUNIPRO} 
proves that the canonical morphism 
$$ 
\Hom_{{\Funct}(\c C_S,\Ind\c D)}(r_{S}(F),G)\R\W
\Hom_{{\Funct}(\c C,\Ind\c D)}(F,GQ)\ .
$$ 
is a bijection for any $G:\c C_S\ra\Ind\c D$. 
\endgraf 
Dually, the functor 
$$l_S:\Funct(\c C,\Pro\c D)\lra\Funct(\c C_S,\Pro\c D)$$ 
is a right adjoint for the canonical functor 
$$ 
Q:\Funct(\c C_S,\Pro\c D)\R\W\Funct(\c C,\Pro\c D)
$$ 
given by the composition with $Q$; 
in particular for any $G:\c C_S\ra\Pro\c D$ we have the bijection 
$$ 
\Hom_{{\Funct}(\c C_S,\Pro\c D)}(G,l_{S}(F))\R\W
\Hom_{{\Funct}(\c C,\Pro\c D)}(GQ,F)\ . 
$$ 
}\par 
\nobreak\hfill$\square$

\Defi{NN}{Deligne localized functors}{
Let $F:\c C\lra\c C'$ a functor, 
and $S$, $S'$ saturated multiplicative systems in
$\c C$, $\c C'$ respectively. 
Then the Deligne right localized functor of $F$ 
with respect to $S$ and $S'$ is 
the functor 
$r_{S,S'}(F)\bydef\Ind(Q'F)r_S:\c C_S\lra\Ind\c C'_{S'}$: 
$$ 
\diagram{
\c C &\R\W^{r'_S} &\Ind\c C &\R\W^{\Ind(F)} &\Ind\c C' \cr 
\ln{Q}\D\H &\nea\rn{r_S} & & &\D\H\rn{\Ind(Q')} \cr 
\c C_S & &\hw\R{70}_{r_{S,S'}(F)}\hw & & \Ind\c C'_{S'} &.\cr 
}$$ 
We will also use the functor 
$r_{S,S'}(F)Q=\Ind(Q'F)r'_S:\c C\lra\Ind\c C'_{S'}$, 
which will be denoted by $r'_{S,S'}(F)$. 
\endgraf 
Dually, the Deligne left localized functor of $F$ 
with respect to $S$ and $S'$ is 
$l_{S,S'}(F)\bydef\Pro(Q'F)l_S:\c C_S\lra\Pro\c C'_{S'}$: 
$$ 
\diagram{
\c C &\R\W^{l'_S} &\Pro\c C &\R\W^{\Pro(F)} &\Pro\c C' \cr 
\ln{Q}\D\H &\nea\rn{l_S} & & &\D\H\rn{\Pro(Q')} \cr 
\c C_S & &\hw\R{70}_{l_{S,S'}(F)}\hw & & \Pro\c C'_{S'} &.\cr 
}$$ 
We will also use  the functor 
$l_{S,S'}(F)Q=\Pro(Q'F)l'_S:\c C\lra\Pro\c C'_{S'}$, 
which will be denoted by $l'_{S,S'}(F)$. 
} 

\medskip 
{\numfont\NNN. \proc } 
We can summarize the situation in the following diagram: 
$$ 
\def\DQ{\vbox to0pt{\hbox to2pt{\hss$\Bigg\downarrow$\hss}\vss}}
\diagram{ 
       &&\Ind\c C &&\hw\R\WW^{\Ind(F)}\hw     &&\Ind\c C'         \cr
&\ln{r'_S}\nea\!\!\!\!\nea\rn{i}  &\DQ&   &&\nea\rn{i'}      &\DQ      \cr
\c C   &&\hw\R\WW^{\qquad F}\hw         &&\c C'  &&\rn{~\Ind(Q')}  \cr 
\DQ    && &&  \DQ                 \cr
\ln{Q~}   &&\Ind\c C_S &&\hw\R\ww\hw         &&\Ind\c C'_{S'}  \cr
&\ln{i}\nea    &&         &&\nea\rn{i'}                  \cr 
\c C_S &&         &&\c C'_{S'}                \cr
}$$ 
where we have the following commutativities: 
$$ 
\eqalign{
r_S Q &=r'_S \cr 
r_{S,S'}(F) &\bydef\Ind(Q'F)r_S \cr 
r'_{S,S'}(F) &\bydef r'_{S,S'}(F) Q=\Ind(Q'F)r_S Q=\Ind(Q'F)r_S \cr 
\Ind(Q'F)i &=\Ind(Q')i'F=i'Q'F \ .
}$$ 
We also have  the extension of $r_{S,S'}(F)$ to $\Ind\c C_{S}$, 
which will be denoted $\ov r_{S,S'}(F)$, 
and the extension of $r'_{S,S'}(F)$ to $\Ind\c C$, 
which will be denoted $\ov r'_{S,S'}(F)$. 
\endgraf 
Dually, we have the following diagram: 
$$ 
\def\DQ{\vbox to0pt{\hbox to2pt{\hss$\Bigg\downarrow$\hss}\vss}}
\diagram{ 
       &&\Pro\c C &&\hw\R\WW^{\Pro(F)}\hw     &&\Pro\c C'         \cr
&\ln{l'_S}\nea\!\!\!\!\nea\rn{i}  &\DQ&   &&\nea\rn{i'}      &\DQ      \cr
\c C   &&\hw\R\WW^{\qquad F}\hw         &&\c C'  &&\rn{~\Pro(Q')}  \cr 
\DQ    && &&  \DQ                 \cr
\ln{Q~}   &&\Pro\c C_S &&\hw\R\ww\hw         &&\Pro\c C'_{S'}  \cr
&\ln{i}\nea    &&         &&\nea\rn{i'}                  \cr 
\c C_S &&         &&\c C'_{S'}                \cr
}$$ 
where we have the following commutativities: 
$$ 
\eqalign{
l_S Q &=l'_S \cr 
l_{S,S'}(F) &\bydef\Pro(Q'F)l_S \cr 
l'_{S,S'}(F) &\bydef l_{S,S'}(F) Q =\Pro(Q'F)l_S Q=\Pro(Q'F)l_S \cr 
\Pro(Q'F)i &=\Pro(Q')i'F=i'Q'F \ .
}$$ 
We also have the extension of $l_{S,S'}(F)$ to $\Pro\c C_{S}$, 
which will be denoted $\ov l_{S,S'}(F)$, 
and the extension of $l'_{S,S'}(F)$ to $\Pro\c C$, 
which will be denoted $\ov l'_{S,S'}(F)$. 

\smallskip 
{\numfont\NNN. \proc } 
The morphism $\delta_S:i\ra r'_S$ induces a morphism 
$$\delta_{S,S'}(F)\bydef\Ind(Q'F)\bullet\delta_S:
\Ind(Q'F)i\lra\Ind(Q'F)r'_S$$ 
of functors 
$\Ind(Q'F)i=i'Q'F$  %%%:\c C\ra\Ind\c C'_{S'}$ 
and 
$r'_{S,S'}(F):\c C\ra\Ind\c C'_{S'}$. 
We denote by $\ov\delta_{S,S'}(F)$ the extended morphism 
between $\Ind(Q'F)=\ov{i'Q'F}$ 
and  $\ov r'_{S,S'}(F):\Ind\c C\ra\Ind\c C'_{S'}$.
\endgraf 
Dually, 
the morphism $\sigma_S:l'_S\ra i$ induces a morphism 
$$\sigma_{S,S'}(F)\bydef\Pro(Q'F)\bullet\sigma_S:
\Pro(Q'F)l'_S\lra\Pro(Q'F)i$$ 
of functors 
$l'_{S,S'}(F):\c C\ra\Pro\c C'_{S'}$
and 
$\Pro(Q'F)i=i'Q'F$.  %%%:\c C\ra\Pro\c C'_{S'}$ 
We denote by $\ov\sigma_{S,S'}(F)$ the extended morphism 
between $\ov l'_{S,S'}(F):\Pro\c C\ra\Pro\c C'_{S'}$
and 
$\Pro(Q'F)=\ov{i'Q'F}$.

\defi{NNN}{
We say that an object $X$ of $\c C$ is %%%inert for $r_SF$ or 
right inert for $F$ (with respect to $S$ and $S'$) 
if $\delta_{S,S'}(F)(X)$ is an isomorphism. 
In particular, right inert objects for $S$ are right inert for any 
functor (with respect to $S$ and any $S'$).
\endgraf 
Dually, we say that an object $X$ of $\c C$ is %%inert for $l_SF$ or 
left inert for $F$ (with respect to $S$ and $S'$) 
if $\sigma_{S,S'}(F)(X)$ is an isomorphism. 
In particular, left inert objects for $S$ are left inert for any 
functor (with respect to $S$ and any $S'$).} 
\rif{DEFINERTOBJFUN}

In the next paragraph we will discuss  the property of 
being localizable for $F$ with respect to $S$ and $S'$. 

\Theo{NNN}{Universal property of Deligne localized functors}{ 
For any $G:\c C_S\lra\Ind\c C'_{S'}$ 
the map $\beta$ to $(\beta\bullet Q)\circ\delta_{S,S'}(F)$ 
induces a bijection 
$$ 
\Hom_{{\Funct}(\c C_S,\Ind\c C'_{S'})}(r_{S,S'}(F),G)\R\W
\Hom_{{\Funct}(\c C,\Ind\c C'_{S'})}(i'Q'F,GQ)\ .
$$ 
In particular for any $G:\c C_S\lra\c C'_{S'}$ we have a bijection 
$$ 
\Hom_{{\Funct}(\c C_S,\Ind\c C'_{S'})}(r_{S,S'}(F),i'G)\R\W
\Hom_{{\Funct}(\c C,\c C'_{S'})}(Q'F,GQ) 
$$ 
induced in the same way. 
\endgraf 
Dually, for any $G:\c C_S\lra\Pro\c C'_{S'}$ 
the map $\beta$ to $\sigma_{S,S'}(F)\circ(\beta\bullet Q)$ 
induces a bijection 
$$ 
\Hom_{{\Funct}(\c C_S,\Pro\c C'_{S'})}(G,l_{S,S'}(F))\R\W
\Hom_{{\Funct}(\c C,\Pro\c C'_{S'})}(GQ,i'Q'F)\ .
$$ 
In particular for any $G:\c C_S\lra\c C'_{S'}$ we have a bijection 
$$ 
\Hom_{{\Funct}(\c C_S,\Pro\c C'_{S'})}(i'G,l_{S,S'}(F))\R\W
\Hom_{{\Funct}(\c C,\c C'_{S'})}(GQ,Q'F) 
$$ 
induced in the same way. 
}

{\proc Proof.}
In fact, this is a consequence of \cite{GENLOCALIZINGUNIPRO}. 
\nobreak\hfill$\square$

\medskip
{\numfont\NN. \proc Composition of Deligne localized functors.} 
The universal property allows us to find canonical morphisms 
for the composite of Deligne localized functors; 
let $F':\c C'\lra\c C''$ be another functor and $S''$ a saturated 
multiplicative system in $\c C''$. We then have a canonical 
morphism: 
$$ 
\delta_{S,S',S''}(F',F):r_{S,S''}(F'F) 
\R\WW 
\ov r_{S',S''}(F') r_{S,S'}(F) 
$$ 
given explicitly by 
$$ 
\delta_{S,S',S''}(F',F)
=\Ind(Q''F')\bullet\ov\delta_{S',S''}\bullet\ov r_S(F) \ .
$$ 
In fact the morphism 
$$ 
\ov\delta_{S',S''}(F')\bullet\delta_{S,S'}(F): 
i''Q''F'F \R\W \ov r_{S',S''}(F') r_{S,S'}(F)Q 
$$
factorizes through the canonical morphism 
$\delta_{S,S''}(F'F):i''Q''F'F\lra r_{S,S''}(F'F)Q$. 
\endgraf 
Dually, we have a canonical 
morphism: 
$$ 
\sigma_{S,S',S''}(F',F):l_{S,S''}(F'F) 
\R\WW 
\ov l_{S',S''}(F') l_{S,S'}(F) 
$$ 
given explicitly by 
$$ 
\sigma_{S,S',S''}(F',F)
=\Pro(Q''F')\bullet\ov\sigma_{S',S''}\bullet\ov l_S(F) \ .
$$

The obvious compatibilities induced by  uniqueness give 
the following result. 

\prop{NNN}{ 
The Deligne right (resp. left) localization is 
a normalized lax 2-functor between the 2-category of 
``categories with multiplicative systems'' 
and the 2-category of ``Ind-categories'' (resp. ``Pro-categories'')
sending $(\c C,S)$ to $\Ind(\c C_S)$ (resp. $\Pro(\c C_S)$)
and $F:\c C\ra\c C'$ to the canonical extension 
of $r_{S,S'}(F)$ (resp. $l_{S,S'}(F)$) to 
the $\Ind$-category (resp. $\Pro$-category); 
the canonical morphisms 
$c_{F',F}\bydef\ov\delta_{S,S',S''}(F',F)$ 
(resp. $c_{F',F}\bydef\ov\sigma_{S,S',S''}(F',F)$)
give the constraints of composition.} 

Remark that the 2-categories are (completely) full subcategories of 
the 2-category of ``categories'' (in particular we do not impose 
any condition on the functors, i.e. the 1-morphisms). 
We have to verify the compatibility (associativity) 
of the constraints of composition, 
i.e. the commutativity of the following diagram 
$$ 
\diagram{
\ov r(F''F'F) &\hw\R{50}^{c_{F'',F'F}}\hw &\ov r(F'')\ov r(F'F) \cr 
\ln{c_{F''F',F}}\D\H &&\D\H\rn{\ov r(F'')\bullet c_{F',F}} \cr 
\ov r(F''F')\ov r(F) &\R\W_{c_{F'',F'}\bullet\ov r(F)} 
& \ov r(F'')\ov r(F')\ov r(F) \cr 
}$$
which is clear because the morphisms involved are defined by universal
properties, so that they are unique; 
and the functoriality at the level of 2-morphisms (i.e. morphisms 
of functors), which is a tedious but elementary verification. 
\hfill$\square$

\prop{NN}{
Let $F:\c C\ra\c C'$ and $G:\c C'\ra\c C$ adjoint functors, 
i.e. 
$$ 
\Hom_{\c C'}(FX,X')\isom\Hom_{\c C}(X,GX') 
$$ 
functorially in objects $X$ of $\c C$ and objects $X'$  of $\c C'$; 
then the Deligne localized functors 
$l_{S,S'}(F):\c C_S\lra\Pro\-\c C'_{S'}$ 
and $r_{S',S}(G):\c C'_{S'}\ra\Ind\-\c C_S$ 
are generalized adjoint functors, i.e.,  we have 
$$ 
\Hom_{\Pro(\c C'_{S'})}(l_{S,S'}(F)X,X')\isom
\Hom_{\Ind(\c C_S)}(X,r_{S,S'}(G)X') 
$$ 
functorially in objects $X$  of $\c C$ and objects $X'$  of $\c C'$. 
}\rif{GENLOCADJ}

{\proc Proof.}
The proof is an easy consequence of commutativity of 
inductive limits in the category of sets: 
$$ 
\eqalign{
\Hom_{\Pro(\c C'_{S'})}(l_{S,S'}(F)X,X')&= 
\Hom_{\Pro(\c C'_{S'})}(\uplimpr_{Y{\ra}X}FY,X') \cr &= 
\limin_{Y{\ra}X}\Hom_{\c C'_{S'}}(FY,X') \cr &\isom 
\limin_{Y{\ra}X}\limin_{X'{\ra}Y'}\Hom_{\c C'}(FY,Y') \cr &\isom 
\limin_{X'{\ra}Y'}\limin_{Y{\ra}X}\Hom_{\c C'}(Y,GY') \cr &\isom 
\limin_{X'{\ra}Y'}\Hom_{\c C_S}(X,GY') \cr &\isom
\Hom_{\Ind(\c C_S)}(X,\uplimin_{X'{\ra}Y'}GY') \cr &= 
\Hom_{\Ind(\c C_S)}(X,r_{S,S'}(G)X') \ .
}$$ 
\nobreak\hfill $\square$

\medskip 
{\numfont\NN. \proc Example: the $\Hom$ bifunctor.} 
Let $\c C$ be a category 
and $S$ a saturated multiplicative system in $\c C$. 
Consider the bifunctor 
$$ 
\Hom_{\c C}:\c C^\opp\times\c C\R\W\Set \ ; 
$$ 
its right localized functor is 
$r_S\Hom_{\c C}\bydef\Ind(\Hom_{\c C})\circ r_S$ 
where 
$$ 
\Ind(\Hom_{\c C}): 
\Ind(\c C^\opp\times\c C)\isom
\Pro(\c C)^\opp\times\Ind(\c C)\R\W\Ind\Set 
$$ 
sends $(X_i,Y_j)$ to $\uplimin_{i,j}\Hom_{\c C}(X_i,Y_j)$, 
and 
$$ 
r_S:(\c C^\opp\times\c C)_{S^\opp\times S}\R\W 
\Ind(\c C^\opp\times\c C)\isom
\Pro(\c C)^\opp\times\Ind(\c C) 
$$ 
sends $(X,Y)$ to $(S/X,Y/S)$. 
Therefore we have 
$$ 
r_S\Hom_{\c C}(X,Y)= 
\Ind(\Hom_{\c C})(S/X,Y/S)= 
\uplimin_{S/X,Y/S}\Hom_{\c C}(X',Y') 
$$ 
and in particular we obtain that 
$$ 
\limin\circ r_S\Hom_{\c C}=\Hom_{\c C_S} \ . 
$$ 

On the other hand, the left localization is defined as 
$l_S\Hom_{\c C}\bydef\Pro(\Hom_{\c C})\circ l_S$ 
where 
$$ 
\Pro(\Hom_{\c C}): 
\Pro(\c C^\opp\times\c C)\isom
\Ind(\c C)^\opp\times\Pro(\c C)\R\W\Pro\Set 
$$ 
sends $(X_i,Y_j)$ to $\uplimpr_{i,j}\Hom_{\c C}(X_i,Y_j)$, 
and 
$$ 
l_S:(\c C^\opp\times\c C)_{S^\opp\times S}\R\W 
\Pro(\c C^\opp\times\c C)\isom
\Ind(\c C)^\opp\times\Pro(\c C) 
$$ 
sends $(X,Y)$ to $(X/S,S/Y)$. 
Therefore we have 
$$ 
l_S\Hom_{\c C}(X,Y)= 
\Pro(\Hom_{\c C})(X/S,S/Y)= 
\uplimpr_{X/S,S/Y}\Hom_{\c C}(X',Y') 
$$ 
and in particular we obtain that 
$$ 
\limpr\circ l_S\Hom_{\c C}(X,Y)\isom 
\limpr_{X/S,S/Y}\Hom_{\c C}(X',Y') 
\isom\Hom_{\c C}(\limin_{X/S}X',\limpr_{S/Y}Y') 
\isom\Hom_{\c C}(\limin r_SX,\limpr l_SY) \ . 
$$

%%%%%%%%%%%%%%%%%%%%%%%%%%%%%%%%%%%%%%%%%%%%%%%%%%%%%%%%%%
\section{Grothendieck-Verdier localized functors. }
\rif{SECGROVER}
%%%%%%%%%%%%%%%%%%%%%%%%%%%%%%%%%%%%%%%%%%%%%%%%%%%%%%%%%%

Let $F:\c C\ra\c C'$ be a functor, 
$S$ and $S'$ be right 
(resp. left for the dual assertions) 
quasi-saturated multiplicative systems in 
$\c C$ and $\c C'$, respectively. 

\defi{NN}{
We say that $R_{S,S'}(F)(X)$ exists, 
or $F$ is right localizable on $X$ with respect to $S$ and $S'$, 
or again that $X$ is right localizable for $F$ with respect to $S$ and $S'$,
if $r_{S,S'}(X)$ is an essentially constant object in $\Ind\c C'_{S'}$, 
i.e. if it is isomorphic in $\Ind\c C'_{S'}$ to an object of 
(the image of) $\c C'_{S'}$. 
This means that there exists $R_{S,S'}(F)(X)$ 
(necessarily isomorphic to 
$\limin r_{S,S'}(F)(X)=\limin_{X/S}F(X')$) in $\c C'_{S'}$ 
such that $i'R_{S,S'}(F)(X)\isom r_{S,S'}(F)(X)$. 
\endgraf 
Dually, we say that $L_{S,S'}(F)(X)$ exists, 
or $F$ is left localizable on $X$ with respect to $S$ and $S'$, 
or again that $X$ is left localizable for $F$ with respect to $S$ and $S'$,
if $l_{S,S'}(X)$ is an essentially constant object in $\Pro\c C'_{S'}$, 
i.e., if it is isomorphic in $\Pro\c C'_{S'}$ to an object of 
(the image of) $\c C'_{S'}$. 
This means that there exists $L_{S,S'}(F)(X)$ 
(necessarily isomorphic to 
$\limpr l_{S,S'}(F)(X)=\limpr_{S/X}F(X')$) in $\c C'_{S'}$ 
such that $i'L_{S,S'}(F)(X)\isom l_{S,S'}(F)(X)$. 
}

Notice that if $X$ is right localizable with respect to $S$ 
(see the definition \cite{DEFLOCOBJ}), 
then any $F$ is right localizable on $X$ and 
$r_{S,S'}F(X)\isom F(R_{S}(X))$. 
Moreover, if $X$ is right inert for $F$ 
with respect to $S$ and $S'$ 
(see the definition \cite{DEFINERTOBJFUN}), 
then ($F$ is right localizable on $X$ and)  
$r_{S,S'}F(X)\isom F(X)$. 
\endgraf 
Dually, if $X$ is left localizable with respect to $S$, 
then any $F$ is left localizable on $X$ and 
$l_{S,S'}F(X)\isom F(L_{S}(X))$. 
Moreover, if $X$ is left inert for $F$ with respect to $S$ and $S'$, 
then ($F$ is left localizable on $X$ and)  
$l_{S,S'}F(X)\isom F(X)$.

\defi{NN}{
We say that the Grothendieck-Verdier localized functor 
$R_{S,S'}F$ exists if for any object $X$ of $\c C$, 
$F$ is right localizable on $X$, 
i.e. if and only if we have a diagram 
$$ 
\diagram{
\c C_S &\R\WW^{r_{S,S'}F} &\Ind\c C'_{S'} \cr 
       &\ln{R_{S,S'}F}\sea \qquad\nea\rn{i'} &    \cr 
       &\c C'_{S'}                &     \cr 
}$$ 
which is commutative up to an isomorphism of functors; 
alternatively: 
$r_{S,S'}(F)$ factorizes through $i'$ up to an isomorphism 
$$ 
\rho_{S,S'}(F):r_{S,S'}(F)\R\W^{\isom}i'R_{S,S'}(F) \ .
$$ 
Dually, we say that the Grothendieck-Verdier localized functor 
$L_{S,S'}F$ exists if for any object $X$ of $\c C$, 
$F$ is left localizable on $X$, 
i.e. iff we have a diagram 
$$ 
\diagram{
\c C_S &\R\WW^{l_{S,S'}F} &\Pro\c C'_{S'} \cr 
       &\ln{L_{S,S'}F}\sea \qquad\nea\rn{i'} &    \cr 
       &\c C'_{S'}                &     \cr 
}$$ 
which is commutative up to an isomorphism of functors; 
alternatively: 
$l_{S,S'}(F)$ factorizes through $i'$ up to an isomorphism 
$$ 
\lambda_{S,S'}(F):l_{S,S'}(F)\R\W^{\isom}i'L_{S,S'}(F) \ .
$$ 
}

We now prove that the 
Grothendieck-Verdier localized functors 
are characterized by the usual universal properties. 

\prop{NN}{ 
The right Grothendieck-Verdier localized functor 
$R_{S,S'}(F):\c C_S\ra\c C'_{S'}$ 
is defined by the following universal property: 
there exists a canonical morphism 
$$\Delta_{S,S'}(F):Q'F\lra R_{S,S'}(F)Q $$ 
such that 
for any $G:\c C_S\ra\c C'_{S'}$ the map 
$$ 
\Hom_{\Funct(\c C_S,\c C'_{S'})}(R_{S,S'}(F),G)\R\W
\Hom_{\Funct(\c C,\c C'_{S'})}(Q'F,GQ)
$$ 
sending $\alpha$ to $(\alpha\bullet Q)\circ\Delta_{S,S'}(F)$ 
is a bijection. 
\endgraf 
Dually, the left Grotendieck-Verdier localized functor 
$L_{S,S'}(F):\c C_S\ra\c C'_{S'}$ 
is defined by the following universal property: 
there exists a canonical morphism 
$$\Sigma_{S,S'}(F):L_{S,S'}(F)Q \lra Q'F $$ 
such that 
for any $G:\c C_S\ra\c C'_{S'}$ the map 
$$ 
\Hom_{\Funct(\c C_S,\c C'_{S'})}(G,L_{S,S'}(F))\R\W
\Hom_{\Funct(\c C,\c C'_{S'})}(GQ,Q'F)
$$ 
sending $\alpha$ to $\Sigma_{S,S'}(F)\circ(\alpha\bullet Q)$ 
is a bijection. 
}

{\proc Proof.}
Suppose that $R_{S,S'}(F)$ exists; 
then there exists an isomorphism 
$\rho_{S,S'}(F):r_{S,S'}(F)\lra^{\isom}i'R_{S,S'}(F)$ 
and the composition 
$$ 
i'Q'F\R\W^{\delta_{S,S'}(F)}
r_{S,S'}(F) Q\R\W^{\rho_{S,S'}(F)\bullet Q} 
i'R_{S,S'}(F)Q 
$$ 
gives a morphism necessarily of the form $i'\bullet\Delta_{S,S'}(F)$ 
with $\Delta_{S,S'}(F):Q'F\lra R_{S,S'}(F)Q $. 
This morphism induces the following bijection: 
$$ 
\eqalign{
\Hom_{\Funct(\c C,\c C'_{S'})}(Q'F,GQ)& \isom 
\Hom_{\Funct(\c C_S,\Ind(\c C'_{S'}))}(r_{S,S'}(F),i'G) \cr 
&\isom \Hom_{\Funct(\c C,\Ind(\c C'_{S'}))}(r_{S,S'}(F)Q,i'GQ) \cr 
&\isom \Hom_{\Funct(\c C,\Ind(\c C'_{S'}))}(i'R_{S,S'}(F)Q,i'GQ) \cr 
&\isom \Hom_{\Funct(\c C,\c C'_{S'})}(R_{S,S'}(F)Q,GQ) \cr 
&\isom \Hom_{\Funct(\c C_S,\c C'_{S'})}(R_{S,S'}(F),G)} 
$$ 
for any $G:\c C_S\lra\c C'_{S'}$,  
which proves the universal property. 
Vice-versa, suppose that there exists a functor 
$R_{S,S'}(F):\c C_S\lra\c C'_{S'}$ endowed with a morphism 
$\Delta_{S,S'}(F):Q'F\lra R_{S,S'}(F)Q $ 
with the stated universal property; 
then by the universal property of $r_{S,S'}(F)$, 
using $G=R_{S,S'}(F)$ we find a canonical morphism 
$$ 
\rho_{S,S'}(F):r_{S,S'}(F)\R\W i'R_{S,S'}(F) 
$$ 
corresponding to $\Delta_{S,S'}(F):Q'F\lra R_{S,S'}(F)Q $ by 
$(\rho_{S,S'}(F)\bullet Q)\circ\delta_{S,S'}(F)
=i'\bullet\Delta_{S,S'}(F)$; 
we have to prove that it is an isomorphism. 
We see that the composite bijection 
$$ 
\eqalign{
\Hom_{\Funct(\c C_S,\Ind(\c C'_{S'}))}(r_{S,S'}(F),i'G)&\isom 
\Hom_{\Funct(\c C,\c C'_{S'})}(Q'F,GQ) \cr 
&\isom \Hom_{\Funct(\c C_S,\c C'_{S'})}(R_{S,S'}(F),G) \cr 
&\isom \Hom_{\Funct(\c C_S,\Ind(\c C'_{S'}))}(i'R_{S,S'}(F),i'G) } 
$$ 
is induced by the composition with $\rho_{S,S'}(F)$ for any $G$; 
this implies that $\rho_{S,S'}(F)$ is an isomorphism, 
so that $R_{S,S'}(F)$ is the Grothendieck-Verdier localized functor of
$F$. 
\nobreak\hfill$\square$ 

\smallskip 
The universal property of the right (resp. left) 
Grothendieck-Verdier localized functor
gives, as in the case of Deligne functors,  
a canonical composition morphism 
$$ 
\Delta_{S,S',S''}(F',F):R_{S,S''}(F'F) \R\W R_{S',S''}(F')R_{S,S'}(F) 
$$ 
(resp. 
$$ 
\Sigma_{S,S',S''}(F',F):R_{S,S''}(F'F) \R\W R_{S',S''}(F')R_{S,S'}(F) \ ),
$$ 
if the terms exist, satisfying the usual associative property. 
Therefore we have the following result. 

\prop{NN}{
The right (resp. left) Grothendieck-Verdier localization gives a 
partially defined, normalized lax 2-functor from the 2-category 
of ``categories with multiplicative systems'' 
to the 2-category of ``categories'' 
sending $(\c C,S)$ to $\c C_S$ 
and $F:\c C\ra\c C'$ to 
$R_{S,S'}(F):\c C_S\lra\c C'_{S'}$ 
(resp. to $L_{S,S'}(F):\c C_S\lra\c C'_{S'}$). 
} 

We are interested in cases in which the canonical morphisms 
$\delta_{S,S',S''}(F',F)$ and $\Delta_{S,S',S''}(F',F)$ 
(resp. $\sigma_{S,S',S''}(F',F)$ and $\Sigma_{S,S',S''}(F',F)$)
for the
compositions are isomorphisms; 
in that case we say, with a slight
abuse of language, that the 2-functors are strict on $F$ and $F'$. 

\prop{NN}{
Suppose that $R_{S,S'}(F)$ and $R_{S',S''}(F')$ exist 
(i.e., that there exist the isomorphisms 
$\rho_{S,S'}(F)$ and $\rho_{S',S''}(F')$); 
then $c_{F',F}$ is an isomorphism if and only if  
$R_{S,S''}(F'F)$ exists 
(i.e. there exists the isomorphism $\rho_{S,S''}(F'F)$) 
and $\Delta_{S,S',S''}(F',F)$ is an isomorphism. 
In particular where $R$ and $r$ are defined as 2-functors, 
one is strict if and only if the other is. 
\endgraf 
Dually, suppose that $L_{S,S'}(F)$ and $L_{S',S''}(F')$ exist 
(i.e., that there exist the isomorphisms 
$\gamma_{S,S'}(F)$ and $\gamma_{S',S''}(F')$); 
then $c_{F',F}$ is an isomorphism if and only if  
$L_{S,S''}(F'F)$ exists 
(i.e. there exist the isomorphism $\gamma_{S,S''}(F'F)$) 
and $\Sigma_{S,S',S''}(F',F)$ is an isomorphism. 
In particular when $L$ and $l$ are defined, 
one is strict if and only if the other is. 
} 

{\proc Proof.}
This is an easy consequence of the following commutative diagram 
$$ 
\diagram{ 
r_{S,S'}(F'F) &\hw\R{80}^{\rho_{S,S''}(F'F)}\hw &i''R_{S,S''}(F'F) \cr 
\ln{\delta_{S,S',S''}(F',F)}\D\H & &\D\H\rn{i''\bullet\Delta_{S,S',S''}(F',F)} \cr 
\ov r_{S',S''}(F') r_{S,S'}(F) &\R\W_{\ov\rho_{S',S''}(F')\rho_{S,S'}(F)} 
&i''R_{S',S''}(F')R_{S,S'}(F) \cr 
}$$ 
\hfill$\square$

\smallskip 
{\numfont\NN. \proc Existence conditions.}\rif{REMEXICON}
Let $S$ be a right (resp. left) 
saturated multiplicative system in $\c C$, 
and let $\c B$ be a full subcategory of $\c C$; 
consider the following conditions: 
\item{$(i)$} 
$\c B$ is right sufficient for $S$, that is 
for any $X$ in $\c C$ there exists $s:X\ra X'$ in $S$ 
with $X'$ an object of $\c B$ 
(resp. $\c B$ is left sufficient for $S$, that is 
for any $X$ in $\c C$ there exists $s:X'\ra X$ in $S$ 
with $X'$ an object of $\c B$); 
\item{$(ii)$} 
for any $X'$ in $\c B$, any $s:X'\ra X$ in $S$ is an isomorphism 
(resp. any $s:X\ra X'$ in $S$ is an isomorphism); 
\item{$(iii)$}  
any morphism in $S$ between objects of $\c B$ is an isomorphism. 
\endgraf\noindent 
We have that $(ii)$ implies $(iii)$, and 
$(iii)$ implies that any object of $\c B$ 
is right (resp. inert) inert for $S$. 
\endgraf\noindent  
The condition $(i)$ implies that 
$$ 
r_{S}(X)\isom\uplimin_{X/S\cap\c B}X' 
\qquad \hbox{(resp. }
l_{S}(X)\isom\uplimpr_{S\cap\c B/X}X' 
\hbox{ )}
$$ 
where $X/S\cap\c B$ (resp. $S\cap\c B/X$) 
is the full subcategory of $X/S$ (resp. $S/X$) whose 
objects are of the form $X\ra X'$ (resp. $X'\ra X$) 
with $X'$ object of $\c B$. 
Moreover, if the condition $(iii)$ holds, then 
any object of $\c C$ is right (resp. left) localizable with respect to $S$, 
and $r_{S}(X)\isom X'$ (resp. $l_{S}(X)\isom X'$) 
for any $X\ra X'$ (resp. $X'\ra X$) 
in $S$ with $X'\in\c B$;  
finally, the right (resp. left) 
inert objects for $S$ are exactly the objects of $\c C$ 
which are isomorphic to some object of $\c B$. 

{\numfont\NNN. \proc } 
In particular if the conditions $(i)$ and $(iii)$ hold, 
then any functor $F:\c C\ra\c C'$ 
is right (resp. left) derivable with respect to 
$S$ and any $S'$. 
In fact $R_{S,S'}(F)X\isom F(X')$ 
(resp. $L_{S,S'}(F)X\isom F(X')$) 
for any $X'$ as before. 

{\numfont\NNN. \proc }
As usual, we may specify the previous conditions 
$(ii)$ and $(iii)$ for a given functor $F:\c C\ra\c C'$ 
and right (resp. left) quasi-saturated multiplicative systems 
$S$ and $S'$ of $\c C$ and $\c C'$, respectively, as follows: 
\item{$(ii(F))$} 
for any $X'$ in $\c B$, and any $s:X'\ra X$ in $S$  
(resp. any $s:X\ra X'$ in $S$), the image $F(s)$ is in $S'$; 
\item{$(iii)$}  
any morphism in $S$ between objects of $\c B$ 
has image (by $F$) in $S'$. 
Clearly, $(ii)$ implies $(iii)$, and 
$(iii)$ implies that any object of $\c B$ 
is right (resp. inert) inert for $F$ with respect to $S$. 
\endgraf\noindent  
Moreover, if $\c B$ is right (resp. left) sufficient for $S$, then 
any object of $\c C$ is right (resp. left) 
localizable for $F$ with respect to $S$ and $S'$, 
and $R_{S,S'}(F)(X)\isom F(X')$ (resp. $L_{S,S'}(F)(X)\isom F(X')$) 
for any $X\ra X'$ (resp. $X'\ra X$) 
in $S$ with $X'\in\c B$.  
% finally, the right (resp. left) 
% inert objects for $F$ with respect to $S$ and $S'$ 
% are exactly the object of $\c C$ 
% which are isomorphic to some object of $\c B$. 

\prop{NN}{
Let $F:\c C\ra\c C'$ and $G:\c C'\ra\c C$ be adjoint functors, 
i.e. 
$$ 
\Hom_{\c C'}(FX,X')\isom\Hom_{\c C}(X,GX') 
$$ 
functorially in objects $X$  of $\c C$ and objects $X'$  of $\c C'$; 
then the Grothendieck-Verdier localized functors  
$L_{S,S'}(F):\c C_S\ra\c C'_{S'}$ and 
$R_{S',S}(G):\c C'_{S'}\ra\c C_S$, 
if they exist,  
are adjoint functors, i.e. we have 
$$ 
\Hom_{\c C'_{S'}}(L_{S,S'}(F)X,X')\isom\Hom_{\c C_S}(X,R_{S',S}(G)X') 
$$ 
functorially in objects $X$ of $\c C$ and  $X'$  of $\c C'$. 
}\rif{LOCADJ}

The proof is an easy consequence of the analogous statement \cite{GENLOCADJ}
for the Deligne localized functors. 
\hfill$\square$

\smallskip

{\numfont\NN. \proc }
As a consequence of the universal properties, 
or of the discussion of Deligne localized functors, 
we remark that the functors 
$R_S$ and $L_S$, if they exist, 
are respectively the left and right adjoint of the canonical functor 
$\Funct(\c C_S,\c D)\ra\Funct(\c C,\c D)$ 
given by the composition with $Q:\c C\ra\c C_S$.

%%%%%%%%%%%%%%%%%%%%%%%%%%%%%%%%%%%%
\section{The case of Triangulated Categories.}
\rif{SECTRICAT}
%%%%%%%%%%%%%%%%%%%%%%%%%%%%%%%%%%%%%%%%%%%%%%%%%%%%%%%%%%

{\numfont\NN. \proc }
If $\c T$ is a triangulated category, 
with translation functor $T$ and class of distinguished triangles 
$\c Tr$, then the categories 
$\Ind(\c T)$ and $\Pro(\c T)$ admit translation functors 
$\Ind(T)$ and $\Pro(T)$ resp., 
and also a special subcategory of the category of (its) triangles, 
give by the essential image of the categories $\Ind(\c Tr)$ 
and $\Pro(\c Tr)$ resp. in the category of triangles of 
$\Ind(T)$ and $\Pro(T)$ respectively. 

But these notions do not give triangulated structures to 
$\Ind(\c T)$ and $\Pro(\c T)$, because two axioms can be violated: 
the completion of morphisms of triangles, and the octahedral axiom. 

\medskip

\lemm{NNN}{
Let $\Delta$ be a distinguished triangle in $\c T$; 
define $\Delta/S$ to be the category whose objects are the 
morphisms of distinguished triangles $\Delta\ra\Delta'$ 
with morphisms in $S$. Then 
$\Delta/S$ is a (right) filtering category 
and we have a canonical isomorphism 
$r_S\Delta\isom\uplimin_{\Delta/S}\Delta'$, 
so that the image by $r_S$ of a distinguished triangle $\Delta$ 
is a triangle of $\Ind(\c T)$ which is an inductive system of 
distinguished triangle of $\c T$; in particular we have a functor 
$$ 
r_S:\c Tr \R\W \Ind(\c Tr)\subseteq \c Tr(\Ind\c T) \ .
$$ 
Dually, 
define $S/\Delta$ to be the category whose objects are the 
morphisms of distinguished triangles $\Delta'\ra\Delta$ 
with morphisms in $S$. 
Then $S/\Delta$ is a left filtering category and  
$l_S\Delta\isom\uplimpr_{S/\Delta}\Delta'$, 
so that the image by $l_S$ of a distinguished triangle $\Delta$ 
is a triangle of $\Pro(\c T)$ which is a projective system of 
distinguished triangles of $\c T$; 
in particular we have a functor 
$$ 
l_S:\c Tr \R\W \Pro(\c Tr)\subseteq \c Tr(\Pro\c T) \ .
$$ 
} 

{\proc Proof.}
The non trivial fact is the property of the category $\Delta/S$, 
which depends on the properties of a saturated multiplicative system, 
especially the completion of squares and the right/left equalization 
of a pair of morphisms. 
Let $s':\Delta\ra\Delta'$ and $s'':\Delta\ra\Delta''$ 
two objects of $\Delta/S$; 
if $\Delta=(X_1\ra X_2\ra X_3\ra TX_1)$, 
$\Delta'=(X'_1\ra X'_2\ra X'_3\ra TX'_1)$ and so on, 
we can complete the commutative diagram 
$$ 
\def\normalbaselines{\baselineskip3pt\lineskip1pt\lineskiplimit1pt}%
\mmatrix{
&&X_2 &\R\ww^{s''}\hw & &X''_2 \cr 
&\ln{f}\nea & &&\nea\rn{f''} \cr 
X_1 &\R\ww\hw &\D\H &X''_1 \cr 
\ln{s'}\D\H & &X'_2 \cr 
&\nea\rn{f'} \cr 
X'_1 \cr 
}
\qquad\hbox{to}\qquad 
\mmatrix{
&&X_2 &\R\ww^{s''}\hw & &X''_2 \cr 
&\ln{f}\nea & &&\ln{f''}\nea \cr 
X_1 &\R\ww\hw &\D\H &X''_1 & &\D\H\rn{u''} \cr 
%% \ln{s'}\D\H & &X'_2 &\R\ww\hw &\D\H &A\cr 
\ln{s'}\D\H & &X'_2 &\D\H &\hw\R\ww &B\cr 
&\nea\rn{f'} \cr 
X'_1  &\R\ww_{u'}\hw & &A\cr 
}
$$ 
where the morphisms named $u$ are in $S$; 
then we can complete the lower side of the cube to a commutative square 
$$ 
\diagram{
& &X'_2 &\R\w^{u'} &B &\R\w^u &C \cr 
&\ln{f'}\nea &&&&\nea\rn{g} \cr 
X'_1 & &\hw\R\WW_{u'}\hw & &A\cr 
}$$ 
with again $u$ in $S$; 
eventually after composition with a morphism $v:C\ra D$ in $S$, 
we obtain a commutative cube 
$$
\def\normalbaselines{\baselineskip1pt\lineskip1pt\lineskiplimit1pt}%
\mmatrix{
&&X_2 &\R\ww^{s''}\hw & &X''_2 \cr 
&\ln{f}\nea & &&\ln{f''}\nea \cr 
X_1 &\R\ww\hw &\D\H &X''_1 & &\D\H\rn{u''} \cr 
%% \ln{s'}\D\H & &X'_2 &\R\ww\hw &\D\H &A\cr 
\ln{s'}\D\H & &X'_2 &\D\H &\hw\R\ww &D\cr 
&\nea\rn{f'} &&&\nea\rn{h} \cr 
X'_1  &\R\ww_{u'}\hw & &A\cr 
}
$$ 
where horizontal and vertical morphisms are in $S$. 
Completing the distinguished triangle on $h$ 
as $\Gamma=(A\ra D\ra E\ra TA)$, 
and the morphisms of the cube to morphisms of triangles, 
we find a diagram of distinguished triangles 
$$ 
\diagram{
\Delta &\R\W^{s''} &\Delta'' \cr 
\ln{s'}\D\h   &     &\D\h\rn{u''}     \cr 
\Delta'&\R\W_{u'} &\Gamma   \cr
}$$ 
in which however the third square 
$$ 
\diagram{
X_3 &\R\W &X''_3 \cr 
\D\h   &     &\D\h     \cr 
X'_3&\R\W &E   \cr
}$$ 
need not be commutative; 
we can choose two morphisms $v',v'':E\ra X'''_3$ in $S$ such that 
$v''u''s''=v'u's'$ and the two compositions 
$D\ra E\ra X'''_3$ coincide (because in any case 
they coincide after composition 
on the right with a morphism of $S$). 
Finally we can take $X'''_2=D\ra X'''_3$, complete the 
distinguished triangle $\Delta'''$ on this morphism, and complete the 
morphisms of distinguished triangles $\Gamma\ra\Delta'''$ 
to obtain a commutative square  
$$ 
\diagram{
\Delta &\R\W^{s''} &\Delta'' \cr 
\ln{s'}\D\h   &     &\D\h\rn{u''}     \cr 
\Delta'&\R\W_{u'} &\Delta'''   \cr
}$$ 
i.e. an object $\Delta\ra\Delta'''$ of $\Delta/S$ with two morphisms 
$\Delta'\ra\Delta'''$ and $\Delta''\ra\Delta'''$. 
\endgraf 
A similar argument, starting with two morphisms 
$f,g:\Delta'\ra\Delta''$ in $\Delta/S$, shows that 
there exists $\Delta\ra\Delta'''$ in $\Delta/S$ with a morphism 
$\Delta''\ra\Delta'''$ of $\Delta/S$ 
which equalizes the two given maps. 
\hfill$\square$ 

\medskip

\prop{NN}{
Let $\Delta=(X_1\ra X_2\ra X_3\ra TX_1)$ be a distinguished triangle in 
$\c T$; 
\item{$(o)$}
if $X_1$ and $X_3$ are right (resp. left) inert or 
localizable with respect to $S$, then so too is $X_2$. 
\endgraf\noindent 
Let $F:\c T\lra\c T'$ be a triangulated functor of triangulated
categories endowed with null systems $\c N$, $\c N'$. 
Then the functor $r_{S,S'}(F)$ (resp.  $l_{S,S'}(F)$) sends 
$\c Tr(\c T/\c N)$ in $\Ind(\c Tr(\c T'/\c N'))$ 
(resp. $\Pro(\c Tr(\c T'/\c N'))$) and 
\item{$(i)$} 
if $X_1$ and $X_3$ are right (resp. left) inert for $F$ 
with respect to $S$ and $S'$, then so too is $X_2$. 
\item{$(ii)$} 
if $R_{S,S'}(F)$ (resp. $L_{S,S'}(F)$) is defined in $X_1$ and $X_3$, 
then it is defined also in $X_2$ and 
$$R_{S,S'}(F)\Delta\bydef 
(R_{S,S'}(F)X_1\ra R_{S,S'}(F)X_2\ra R_{S,S'}(F)X_3\ra R_{S,S'}(F)TX_1)$$ 
(resp. 
$$L_{S,S'}(F)\Delta\bydef 
(L_{S,S'}(F)X_1\ra L_{S,S'}(F)X_2\ra L_{S,S'}(F)X_3\ra L_{S,S'}(F)TX_1)
\quad\hbox{ )}
$$
is a distinguished triangle of $\c T'/\c N'$. 
}

{\proc Proof.}
The arguments are a slight modification of Deligne [SGA4,{\it XVII},1.2.2]. 
For $(i)$ consider the canonical morphism 
$F(\Delta)\ra rF(\Delta)$; the first is a distinguished triangle, 
and the second is an ind-object of distinguished triangles. 
Moreover the first and third morphisms are iso. 
Therefore for any $W$ we have a morphism of long exact sequences 
$$ 
\Hom_{\c T}(W,F(\Delta))\R\W\Hom_{\c T}(W,rF(\Delta))= 
\limin_{\Delta/S}\Hom_{\c T}(W,F(\Delta')) 
$$ 
in which we can apply the five lemma. 
For $(ii)$ we consider the distinguished triangle constructed from the 
morphism $RF(X_3)\ra RF(TX_1)$ and the morphisms 
$$ 
\diagram{ 
RF(X_1)       &\R\W &W &\R\W &RF(X_3)        \cr 
\ln{\isom}\D\h &     &  &    &\D\h\rn{\isom} \cr 
rF(X_1)       &\R\W &rF(X_2) &\R\W &rF(X_3)        \cr 
}$$ 
where the second line is again an ind-object of distinguished
triangles; so we can extend the diagram to a morphism of 
triangles, which is necesserily an isomorphism. 
The point $(o)$ can be deduced from $(i)$ and $(ii)$ 
using $F=\id_{\c T}$. 
\hfill$\square$ 

\medskip

\prop{NN}{
Let $F:\c T\lra\c A$ be a cohomological functor between a triangulated
category $\c T$  endowed with a null systems $\c N$,
and an abelian category $\c A$. 
Then the functor $r_{S}(F):\c T/\c N\lra\Ind(\c A)$ 
(resp.  $l_{S}(F):\c T/\c N\lra\Pro(\c A)$) 
is a cohomological functor. 
}

{\proc Proof.} 
Elementary. 
\hfill$\square$ 

\medskip 

%%%%%%%%%%%%%%%%%%%%%%%%%%%%%%%%%%%%
\section{The case of Derived Categories.}
\rif{SECDERCAT}
%%%%%%%%%%%%%%%%%%%%%%%%%%%%%%%%%%%%%%%%%%%%%%%%%%%%%%%%%%
\def\qis{{\rm qis}}

{\numfont\NN. \proc }
Any additive functor $F:\c A\ra\c A'$ 
of abelian categories extends to 
an additive functor 
$F=\bC^\ast(F):\bC^\ast(\c A)\ra\bC^\ast(\c A')$ 
and to a triangulated functor 
$F=\bK^\ast(F):\bK^\ast(\c A)\ra\bK^\ast(\c A')$ 
of categories of homotopic complexes. 
Localizing with respect to the systems of quasi-isomorphisms 
of $\bC^\ast(\c A)$ and $\bC^\ast(\c A')$, 
we have the right Deligne derived functor 
$$ 
r(F)=r(\bK^\ast(F)):\bD^\ast(\c A)\R\W\Ind(\bD^\ast(\c A')) 
$$ 
defined by the composition 
$$ 
\bD^\ast(\c A)\R\W^r 
\Ind(\bK^\ast(\c A))\R\W^{\Ind(\bK^\ast(F))} 
\Ind(\bK^\ast(\c A'))\R\W^{\Ind(Q')} 
\Ind(\bD^\ast(\c A'))
$$ 
having the universal property 
$$ 
\Hom_{{\Funct}(\bD^\ast(\c A),\Ind\bD^\ast(\c A'))}(r(F),G)
\R\W^{\isom}
\Hom_{{\Funct}(\bK^\ast(\c A),\Ind\bD^\ast(\c A'))}(i'Q'F,GQ)\ ,
$$ 
induced by the canonical morphism 
$\delta(F):i'Q'F\ra r(F)Q$. 
\endgraf 
Dually we have the left Deligne derived functor 
$$ 
l(F)=l(\bK^\ast(F)):\bD^\ast(\c A)\R\W\Pro(\bD^\ast(\c A')) 
$$ 
defined by the composition 
$$ 
\bD^\ast(\c A)\R\W^l 
\Pro(\bK^\ast(\c A))\R\W^{\Pro(\bK^\ast(F))} 
\Pro(\bK^\ast(\c A'))\R\W^{\Pro(Q')} 
\Pro(\bD^\ast(\c A'))
$$ 
having the universal property 
$$ 
\Hom_{{\Funct}(\bD^\ast(\c A),\Pro\bD^\ast(\c A'))}(G,l(F))
\R\W^{\isom}
\Hom_{{\Funct}(\bK^\ast(\c A),\Pro\bD^\ast(\c A'))}(GQ,i'Q'F)\ ,
$$ 
induced by the canonical morphism 
$\sigma(F):l(F)Q\ra i'Q'F$. 

\smallskip
{\numfont\NNN. \proc }
We study the condition for a complex to be 
right (resp. left) inert for the 
multiplicative system of quasi-isomorphisms. 
A sufficient condition is 
to be homotopically equivalent to a complex whose terms 
are injective (resp. projective) objects of 
the abelian category $\c A$. 
The condition is also necessary if the category 
$\c A$ has enough injective (resp. projective) objects. 
In fact the subcategory of these objects satisfies the 
conditions of \cite{REMEXICON} 
(see for example [RD,I.4.5]).

\smallskip
{\numfont\NNN. \proc }
We have also the notion of right and left 
Grothendieck-Verdier derived functors which are defined if 
the  corresponding Deligne derived functors take their image in the 
essentially constant objects. 

\smallskip

{\numfont\NN. \proc } 
As an example we consider the $\Hom^\point$ functor: 
$$ 
\Hom^\point:\bC^b(\c A)^\opp\times\bC^b(\c A)\R\W\bC^b(\Ab)
$$ 
for which we have $H^0(\Hom^\point(X,Y))=\Hom_{\bK^b(\c A)}(X,Y)$. 
The right Deligne derived functor is defined as 
$$ 
r\Hom^\point\bydef\Ind(\Hom^\point)r: 
\bD^b(\c A)^\opp\times\bD^b(\c A)\R\W\Ind(\bD^b(\Ab))
$$ 
sending $(X,Y)$ to $\uplimin_{\qis/X,Y/\qis}\Hom^\point(X',Y')$. 
Composition with the inductive limit commutes with passage to 
cohomology and we have 
$$ 
H^0(\limin_{\qis/X,Y/\qis}\Hom^\point(X',Y'))\isom 
\limin_{\qis/X,Y/\qis}\Hom_{\bK^b(\c A)}(X',Y')= 
\Hom_{\bD^b(\c A)}(X,Y) \ .
$$ 
By contrast, the left Deligne derived functor is defined as 
$$ 
l\Hom^\point\bydef\Pro(\Hom^\point)l: 
\bD^b(\c A)^\opp\times\bD^b(\c A)\R\W\Pro(\bD^b(\Ab))
$$ 
sending $(X,Y)$ to $\uplimpr_{X/\qis,\qis/Y}\Hom^\point(X',Y')$. 
The zero-cohomology gives 
$$ 
H^0(\uplimpr_{\qis/X,Y/\qis}\Hom^\point(X',Y'))\isom 
\uplimpr_{\qis/X,Y/\qis}\Hom_{\bK^b(\c A)}(X',Y') \ .
$$ 
while composition with the projective limit gives 
$$ 
\limpr l\Hom^\point(X,Y)= 
\limpr_{\qis/X,Y/\qis}\Hom^\point(X',Y')
=\Hom^\point(\limin(\qis/X),\limpr(Y/\qis)) \ .
$$

\medskip

{\numfont\NN. 
\proc Comparison between $\bD^\ast(\Ind(\c A))$ and $\Ind(\bD^\ast(\c A))$.} 
Since $\Ind(\c A)$ and $\Pro(\c A)$ are abelian categories if $\c A$ is, 
we may work out the usual derived categories 
$\bD^\ast(\Ind(\c A))$ and $\bD^\ast(\Pro(\c A))$ 
which are triangulated categories, in contrast with 
$\Ind(\bD^\ast(\c A))$ and $\Pro(\bD^\ast(\c A))$ which are not, 
but are involved in the definition of Deligne derived functors. 
We want to make explicit the relations between these categories 
when $\ast=b$ (bounded complexes). 
We recall the following result of Deligne (see [RD, App.,prop.~3]). 

\lemm{NNN}{
The family of functors 
$H^p\bydef\Ind H^p:\Ind(\bD^b(\c A))\ra\Ind\Ab$ 
is conservative, that is a morphism $\phi$ in $\Ind(\bD^b(\c A))$ 
is an isomorphism in $\Ind(\bD^b(\c A))$ if and only if  for every $p$ 
the morphism $H^p(\phi)$ is an isomorphism in $\Ind\Ab$. 
}\rif{LEMCONDEL}

{\numfont\NNN. } 
 The canonical functor 
$\c A\ra\bC^\ast(\c A)$ extends to a functor 
$J:\Ind(\c A)\lra\Ind(\bC^\ast(\c A))$, 
and using the representation \cite{PREMORIND} 
for the morphisms in the $\Ind$-categories, 
for $\ast=b$ we may extend this to a functor 
$\bC^b(\Ind(\c A))\lra\Ind(\bC^b(\c A))$. 
Note that an object in $\bC^b(\Ind(\c A))$ must be represented with 
``commutation relations'' $d^2=0$ (it is a complex), where the zero 
morphisms of ind-objects are represented with all components the zero
morphisms of the category $\c A$; this guarantees that the result 
of the parallelization process is an inductive system of 
complexes. 

Using the Deligne lemma, we may define a functor 
$J:\bD^b(\Ind(\c A))\lra\Ind(\bD^b(\c A))$ 
making  the following diagram commutative 
$$\diagram{
\bC^b(\Ind(\c A)) &\R\W^J &\Ind(\bC^b(\c A)) \cr 
\ln{Q}\D\H              &     &\D\H\rn{\Ind(Q)} \cr 
\bD^b(\Ind(\c A)) &\R\W^J &\Ind(\bD^b(\c A)) \cr 
}$$ 

We remark that the essential image of $\bD^b(\Ind(\c A))$ 
in $\Ind(\bD^b(\c A))$ consists of Ind-objects in the category of 
complexes whose cohomology is {\it uniformly} bounded, 
that is, not only is it an inductive system of bounded cohomology 
complexes, but the bound is uniform on the whole of the system. 
This remark motivates the following definitions. 

\defi{NNN}{
Let $\ind(\bC^b(\c A))$ 
(resp. $\ind(\bK^b(\c A))$, $\ind(\bD^b(\c A))$) be
the full subcategory of $\Ind(\bC^b(\c A))$ 
(resp. $\Ind(\bK^b(\c A))$, $\Ind(\bD^b(\c A))$) 
whose objects are inductive systems of complexes $X_\alpha$ 
such that there exists $N\in\b N$ with 
$X_{\alpha,i}=0$ 
(resp. $X_{\alpha,i}=0$, $H^i(X_\alpha)=0$) 
for all $i\notin[-N,N]$ and all $\alpha$. 
\endgraf 
Dually, 
$\pro(\bC^b(\c A))$, $\pro(\bK^b(\c A))$, $\pro(\bD^b(\c A))$ 
are the full subcategories of 
$\Pro(\bC^b(\c A))$, $\Pro(\bK^b(\c A))$, $\Pro(\bD^b(\c A))$ 
having uniformly limited complexes, or cohomology. 
} 

So we have the following canonical functors: 
$$ 
\diagram{
\bC^b(\Ind(\c A)) &\R\W^J &\ind(\bC^b(\c A)) \cr 
\ln{}\D\H              &     &\D\H\rn{} \cr 
\bK^b(\Ind(\c A)) &\R\W &\ind(\bK^b(\c A)) \cr 
\ln{Q}\D\H              &     &\D\H\rn{\Ind(Q)} \cr 
\bD^b(\Ind(\c A)) &\R\W_J &\ind(\bD^b(\c A)) \cr 
} 
\qquad 
\diagram{
\bC^b(\Pro(\c A)) &\R\W^J &\pro(\bC^b(\c A)) \cr 
\ln{}\D\H              &     &\D\H\rn{} \cr 
\bK^b(\Pro(\c A)) &\R\W &\pro(\bK^b(\c A)) \cr 
\ln{Q}\D\H              &     &\D\H\rn{\Ind(Q)} \cr 
\bD^b(\Pro(\c A)) &\R\W_J &\pro(\bD^b(\c A)) &.\cr 
} 
$$ 
Remark however that the second and third horizontal functors 
do not have, in general, good properties (see [KS3]). 

{\numfont\NNN. } 
Similarly, the canonical functor 
$\c A\ra\Ind\c A$ extends to the complexes 
$\bC^\ast(\c A)\ra\bC^\ast(\Ind\c A)$, 
and trivially we have an extension 
$\ind(\bC^b(\c A))\ra\bC^\ast(\Ind\c A)$ 
since any inductive system of uniformly bounded complexes 
can be rewrittten as a bounded complex of inductive systems. 
In fact the functor $J$ is the Ind-adjoint (resp. the Pro-adjoint) 
of the canonical functor $\bC^\ast(\c A)\ra\bC^\ast(\Ind\c A)$ 
(resp. $\bC^\ast(\c A)\ra\bC^\ast(\Pro\c A)$).

\prop{NNN}{
The categories $\ind(\bC^b(\c A))$ and $\bC^b(\Ind(\c A))$ 
are equivalent, %% as well as 
%% $\ind(\bK^b(\c A))$ and $\bK^b(\Ind(\c A))$; 
and dually 
$\pro(\bC^b(\c A))$ %% ,  $\pro(\bK^b(\c A))$ 
is equivalent to $\bC^b(\Pro(\c A))$ %% , $\bK^b(\Pro(\c A))$. 
}

In fact the two functors previous defined are equivalence 
quasi-inverses  of each other. 
\hfill$\square$

\prop{NNN}{
The natural functor 
$\bD^b(\Ind(\c A))\lra\ind(\bD^b(\c A))$ 
(risp. $\bD^b(\Pro(\c A))\lra\pro(\bD^b(\c A))$) 
is conservative. 
%%  it is fully faithful 
%% if $\c A$ has enough injective objects, 
%% but in general it is not essentially surjective. 
} 

This a consequence of \cite{LEMCONDEL}.
\hfill$\square$

\lemm{NNN}{
Consider the inclusions of $\bK^b(\c A)$ 
in $\ind(\bK^b(\c A))$ and $\bK^b(\Ind\c A)$ 
(resp. in $\pro(\bK^b(\c A))$ and $\bK^b(\Pro\c A)$) 
and of $\bD^b(\c A)$  
in $\ind(\bD^b(\c A))$ and $\bD^b(\Ind\c A)$ 
(resp. in $\pro(\bD^b(\c A))$ and $\bD^b(\Pro\c A)$). 
Then for any $X$ object of $\bC^b(\c A)$ and any 
$Z$ in $\bK^b(\Ind\c A)$ (resp. $\bK^b(\Pro\c A)$) 
we have 
$$ 
\Hom_{\bK^b(\Ind\c A)}(X,Z)\R\w^{\isom} 
\Hom_{\ind(\bK^b(\c A))}(X,JZ) 
\qquad\hbox{(resp.}\quad 
\Hom_{\bK^b(\Pro\c A)}(Z,X)\R\w^{\isom} 
\Hom_{\pro(\bK^b(\c A))}(JZ,X) 
\hbox{\ )}
$$ 
and 
$$ 
\Hom_{\bD^b(\Ind\c A)}(X,Z)\R\w^{\isom} 
\Hom_{\ind(\bD^b(\c A))}(X,JZ) 
\qquad\hbox{(resp.}\quad 
\Hom_{\bD^b(\Pro\c A)}(Z,X)\R\w^{\isom} 
\Hom_{\pro(\bD^b(\c A))}(JZ,X) 
\hbox{\ ).}
$$ 
In particular, the Ind-object $JZ$ as a functor is defined by 
$\Hom_{\bK^b(\Ind\c A)}(\hbox{-},Z)$ in $\ind(\bK^b(\c A))$ 
and $\Hom_{\bD^b(\Ind\c A)}(\hbox{-},Z)$ in $\ind(\bD^b(\c A))$ 
(resp. the Pro-object $JZ$ as a functor is defined by 
$\Hom_{\bK^b(\Pro\c A)}(Z,\hbox{-})$ in $\pro(\bK^b(\c A))$ 
and $\Hom_{\bD^b(\Pro\c A)}(Z,\hbox{-})$ in $\pro(\bD^b(\c A))$). 
}

In fact we may suppose that $Z$ in $\bK^b(\Ind\c A)$ is already 
parallelized, and in that case the boundedness condition gives 
an isomorphism 
$$ 
\Hom_{\bK^b(\Ind\c A)}(X,(Z_i))\isom 
\limin_i\Hom_{\bK^b(\c A)}(X,Z_i) 
$$ 
which is by definition $\Hom_{\ind(\bK^b(\c A))}(X,JZ)$. 
\hfill$\square$

%%%%%%%%%%%%%%%%%%%%%%%%%%%%%%%%%%%%%%%%%%%%%%%%%%%%%%%
%\bigskip
\refer

\biblio{AM} 
       {Artin, M. and Mazur, B.; Etale homotopy.}
       {Lecture Notes in Mathematics, \it Springer-Verlag,}
       {\bf 100} (1986).

\biblio{BBD} 
       {Bernestein, Beilison, Deligne; Faisceaux pervers.}
       {\it Ast\`erisque} 
       {\bf 100} (1982).

\biblio{CT}
       {Saavedra Rivano, Neantro; Cat\'egories Tannakiennes.}
       {Lecture Notes in Mathematics, \it Springer-Verlag,}
       {\bf 265} (1972).

\biblio{D}
       {Deligne, Pierre; \it Cat\'egories tannakiennes.}
       {The Grothendieck Festschrift, Vol. II, 111--195, Progr. Math.,}
       {\bf 87} (1990).

\biblio{DM}
       {Deligne, Pierre and Milne, JS; \it Tannakian categories.}
       {Hodge Cycles, Motives, and Shimura Varieties, Lecture Notes in Math.}
       {\bf 900}, 101-228 (1982).

\biblio{EGA}
       {Grothendieck, Alexander and Dieudonn\'e, Jean; 
        El\'ements de g\'eom\'etrie alg\'ebrique.}   
        {\it Inst. Hautes Etudes Sci. Publ. Math.}, 
        {\bf 4} (1960), {\bf 8} (1961), {\bf 11} (1961), {\bf 17} (1963), 
        {\bf 20} (1964), {\bf 24} (1965), {\bf 28} (1966), {\bf 32} (1967).

\biblio{LC}
       {Grothendieck, Alexander; 
        Local cohomology (notes by Hartshorne, Robin).} 
       {Lecture Notes in Mathematics, \it Springer-Verlag,}
       {\bf 41} (1967).
 
\biblio{KS1} 
       {Kashiwara, Masaki and Schapira, Pierre,}
       {Sheaves on manifolds. \it Springer-Verlag, Berlin.}
       {\bf } (1994).

\biblio{KS2} 
       {Kashiwara, Masaki and Schapira, Pierre;}
       {Ind-sheaves. \it Ast\`erisque}
       {\bf 271} (2001).

\biblio{KS3} 
       {Kashiwara, Masaki and Schapira, Pierre;}
       {\it Microlocal study of ind-sheaves: 
        microsupport and regularity.}
       {\bf } Preprint AG/0108065 (2001).

\biblio{RD} 
       {Hartshorne, Robin; Residues and duality
       (With an appendix by P. Deligne).}
       {Lecture Notes in Mathematics, \it Springer-Verlag,}
       {\bf 20} (1966).

\biblio{SGA1} 
       {Rev\^etements \'etales et groupe fondamental 
        (Dirig\'e par Alexandre Grothendieck).}
       {Lecture Notes in Mathematics, \it Springer-Verlag,}
       {\bf 224} (1971).

\biblio{SGA4} 
       {Th\'eorie des topos et cohomologie \'etale des sch\'emas 
       (Dirig\'e par M. Artin, A. Grothendieck et J. L. Verdier. Avec la
       collaboration de  N. Bourbaki, P. Deligne et B. Saint-Donat).}
       {Lecture Notes in Mathematics, \it Springer-Verlag,}
       {\bf 269} (1972), {\bf 270} (1972), {\bf 305} (1973).

\end